\documentclass[12pt,leqno]{article}
\usepackage{amsmath,amssymb,amsthm,amscd,dsfont}
\numberwithin{equation}{section} \allowdisplaybreaks
\begin{document}
\newtheorem{theorem}{Theorem}[section]
\newtheorem{defin}{Definition}[section]
\newtheorem{prop}{Proposition}[section]
\newtheorem{corol}{Corollary}[section]
\newtheorem{lemma}{Lemma}[section]
\newtheorem{rem}{Remark}[section]
\newtheorem{example}{Example}[section]
\title{Generalized CRF-structures}
\author{{\small by}\vspace{2mm}\\Izu Vaisman}
\date{}
\maketitle
{\def\thefootnote{*}\footnotetext[1]%
{{\it 2000 Mathematics Subject Classification: 53C15} .
}}
\begin{center} \begin{minipage}{12cm}
A{\footnotesize BSTRACT. A generalized F-structure is a complex,
isotropic subbundle $E$ of $T_cM\oplus T^*_cM$
($T_cM=TM\otimes_{\mathds{R}}\mathds{C}$ and the metric is defined
by pairing) such that $E\cap\bar E^{\perp}=0$. If $E$ is also
closed by the Courant bracket, $E$ is a generalized CRF-structure.
We show that a generalized F-structure is equivalent with a
skew-symmetric endomorphism $\Phi$ of $TM\oplus T^*M$ that
satisfies the condition $\Phi^3+\Phi=0$ and we express the
CRF-condition by means of the Courant-Nijenhuis torsion of $\Phi$.
The structures that we consider are generalizations of the
F-structures defined by Yano and of the CR (Cauchy-Riemann)
structures. We construct generalized CRF-structures from: a
classical F-structure, a pair $(\mathcal{V},\sigma)$ where
$\mathcal{V}$ is an integrable subbundle of $TM$ and $\sigma$ is a
$2$-form on $M$, a generalized, normal, almost contact structure
of codimension $h$. We show that a generalized complex structure
on a manifold $\tilde M$ induces generalized CRF-structures into
some submanifolds $M\subseteq\tilde M$. Finally, we consider
compatible, generalized, Riemannian metrics and we define
generalized CRFK-structures that extend the generalized K\"ahler
structures and are equivalent with quadruples
$(\gamma,F_+,F_-,\psi)$, where $(\gamma,F_\pm)$ are classical,
metric CRF-structures, $\psi$ is a $2$-form and some conditions
expressible in terms of the exterior differential $d\psi$ and the
$\gamma$-Levi-Civita covariant derivative $\nabla F_\pm$ hold. If
$d\psi=0$, the conditions reduce to the existence of two partially
K\"ahler reductions of the metric $\gamma$. The paper ends by an
Appendix where we define and characterize generalized Sasakian
structures.}
\end{minipage} \end{center} \vspace{5mm}
\section{Introduction}
This paper belongs to the framework of generalized structures on a
differentiable manifold $M$ (e.g., Hitchin
\cite{H}). These structures are similar to classical
structures but they are defined on the big tangent bundle
$T^{big}M=TM\oplus T^*M$ with the neutral metric
\begin{equation}\label{gFinC}
g((X,\alpha),(Y,\beta))=\frac{1}{2}(\alpha(Y)+\beta(X))
\end{equation} and the Courant bracket \cite{C}
\begin{equation}\label{crosetC} [(X,\alpha),(Y,\beta)] = ([X,Y],
L_X\beta-L_Y\alpha+\frac{1}{2}d(\alpha(Y)-\beta(X))).\end{equation}

In (\ref{gFinC}) and (\ref{crosetC}) the notation uses the
following conventions that will be followed throughout the whole
text: $M$ is an $m$-dimensional manifold and $X,Y,..$ are either
contravariant vectors or vector fields, $\alpha,\beta,...$ are
either covariant vectors or $1$-forms. Furthermore, we will denote
by $\chi^k(M)$ the space of $k$-vector fields, by $\Omega^k(M)$
the space of differential $k$-forms, by $\Gamma$ spaces of global
cross sections of vector bundles, by  $d$ the exterior
differential and by $L$ the Lie derivative. All the manifolds and
mappings are assumed of the $C^\infty$ class.

Most of the work done until now in this framework was on
generalized, complex and K\"ahler structures \cite{Galt} and this
work was motivated by applications to supersymmetry in string theory
\cite{LMTZ}. Other generalized structures that were considered are
the almost product structures \cite{{V-gcm},{Wade}}, the almost
tangent structures \cite{V-gcm} and the almost contact structures
\cite{{PW},{V-stable}}.

A generalized, complex structure can be defined as a complex Dirac
structure, i.e., a maximal $g$-isotropic subbundle of the
complexified big tangent bundle
$T^{big}_cM=T^{big}M\otimes\mathds{C}$ with the space of cross
sections closed by the Courant bracket \cite{Galt}. By dropping
the maximality condition, i.e., by replacing the Dirac structure
by a big-isotropic structure \cite{V06}, we will obtain the notion
of a generalized F-structure, such that any classical F-structure
\cite{Y} produces a generalized F-structure. Furthermore, if the
big-isotropic structure is integrable, i.e., its space of cross
sections is closed by the Courant bracket, we get an integrable,
generalized F-structure.

If the generalized F-structure defined by a classical one is
integrable, the classical F-structure is such that its
$\sqrt{-1}$-eigenbundle (the holomorphic distribution) is a CR
(Cauchy-Riemann) structure \cite{DT}. Accordingly, we propose the
name CRF-structure, which explains the title of this
paper\footnote{The names CR-structure and F-structure are well
established in the mathematical literature. It is interesting to
mention that an equivalent definition of classical F-structures
was given in
\cite{Barros} where these structures are called hor-complex
structures; ``hor" comes from ``horizontal". Correspondingly, we
should use the name of a generalized hor-complex structure. The
``hor-complex" terminology is nice but it was not adopted in the
literature.}.

The addition of a compatible, generalized, Riemannian metric leads
to a notion of generalized CRFK-structure that extends the
generalized K\"ahler structures \cite{Galt}. It turns out that a
generalized CRFK-structure is equivalent with a quadruple
$(\gamma,F_+,F_-,\psi)$ where $(\gamma,F_\pm)$ are classical,
metric CRF-structures, $\psi$ is a $2$-form and some conditions
expressible in terms of the exterior differential $d\psi$ and the
$\gamma$-Levi-Civita covariant derivative $\nabla F_\pm$ hold. If
$d\psi=0$, the conditions reduce to the existence of two partially
K\"ahler reductions of the metric $\gamma$.

The connection with the theory of CR-structures is a motivation of
the present paper. The motivation is enhanced by the fact that
some submanifolds of a generalized, complex manifold may have an
induced, generalized F-structure, like in the case of the well
known CR-submanifolds of a Hermitian manifold \cite{Bej}.\\

We recall that the real big-isotropic structures $E\subseteq
T^{big}M$ were studied in
\cite{V06}. The integrability condition of the structure (closure
by the Courant brackets) is
\begin{equation}\label{Cclosed} [(X,\alpha),(Y,\beta)]\in\Gamma(E),
\hspace{3mm}\forall(X,\alpha),(Y,\beta)\in\Gamma(E)\end{equation}
and the properties of the Courant bracket imply that condition
(\ref{Cclosed}) is equivalent with the condition $[E,E']\subseteq
E'$ where $E'=E^{\perp_g}$ is the $g$-orthogonal bundle of $E$ in
$T^{big}M$. The simplest examples are:
\begin{example}\label{exgraphtheta} {\rm Let $F$ be a subbundle of
$TM$ and $\theta\in\Omega^2(M)$. Then,
\begin{equation}\label{Etheta}
E_\theta=graph(\flat_\theta|_{F})=\{(X,\flat_{\theta}X=i(X)\theta)\,/\,X\in
F\}
\subseteq T^{big}M\end{equation} is a big-isotropic structure on $M$ with
the $g$-orthogonal bundle
\begin{equation}\label{perpEtheta} E'_\theta=\{(Y,\flat_{\theta}
Y+\gamma)\,/\,Y\in TM,\,\gamma\in ann\,S\}.\end{equation} The
structure $E_\theta$ is integrable iff $F$ is a foliation and
$\theta$ satisfies the condition \cite{{Galt},{V06}}
\begin{equation}\label{integrEtheta}
d\theta(X_1,X_2,Y)=0,\hspace{1cm}\forall X_1,X_2\in\Gamma
F,Y\in\chi^1(M).\end{equation}}\end{example}
\begin{example}\label{exgraphP} {\rm
Let $\Sigma$ be a subbundle of $T^*M$ and $P\in\chi^2(M)$. Then
\begin{equation}\label{eqEP} E_P=graph(\sharp_P|_{\Sigma}) =
\{(\sharp_P\sigma=i(\sigma)P,\sigma)\,/\,\sigma\in \Sigma\}\end{equation} is
a big-isotropic structure on $M$ with the $g$-orthogonal bundle
\begin{equation}\label{eqE'P} E'_P=\{(\sharp_P\beta+Y,\beta)\,/\,\beta\in T^*M,
Y\in ann\,\Sigma\}.\end{equation} The structure $E_P$ is
integrable iff $\Sigma\subseteq T^*M$ is closed with respect to
the bracket of $1$-forms
\begin{equation}\label{croset1forme} \{\alpha,\beta\}_P=
L_{\sharp_P\alpha}\beta-L_{\sharp_P\beta}\alpha
-d(P(\alpha,\beta))\end{equation} and the Schouten-Nijenhuis bracket
$[P,P]$ (e.g., \cite{V-carte}) satisfies the condition \cite{V06}
\begin{equation}\label{integrEP}
[P,P](\alpha_1,\alpha_2,\beta)=0, \hspace{1cm}\forall\alpha_1,
\alpha_2\in \Sigma,\beta\in\Omega^1(M).\end{equation}}\end{example}

We end the Introduction by a few words about the content of the
paper. In Section 2 we define the generalized CRF-structures as
complex, big-isotropic structures, prove the equivalence with a
skew-symmetric endomorphism $\Phi\in End(T^{big})M$ such that
$\Phi^3+\Phi=0$, which may also be presented as a triple of
classical tensor fields $(A\in\Gamma(TM\otimes
T^*M),\sigma\in\Omega^2(M),\pi\in\chi^2(M))$, and express the
integrability condition in terms of $\Phi$. Finally, we consider
generalized F-structures induced on some submanifolds of a
generalized complex manifold and show their integrability. In
Section 3 we discuss particular classes of CRF-structures. We show
that a classical F-structure produces a generalized F-structure
and get the corresponding integrability conditions, which define a
seemingly new class of Yano's F-structures. Then, we associate a
generalized F-structure with a pair $(\mathcal{V}\subseteq
TM,\theta\in\Omega^2(M))$, where $\theta$ is non degenerate on
$\mathcal{V}$, and with a pair $(\Sigma\subseteq
T^*M,\pi\in\Omega^2(M))$, where $\pi$ is non degenerate on
$\Sigma$. In Section 4, we present the generalized Riemannian
metrics following \cite{Galt} and extend the notion of a
generalized K\"ahler manifold to a notion of CRFK-manifold by
replacing the generalized complex structure by a generalized
CRF-structure. Then we obtain expressions of the CRFK condition by
means of the corresponding classical objects and show that
Riemannian manifolds $(M,\gamma)$ where the metric $\gamma$ has
two partially K\"ahler reductions are CRFK-manifolds.
\section{The basics of generalized CRF-structures}
The definition of big-isotropic structures, integrability and
other basic properties may be complexified, i.e., transferred to
the complexified tangent bundle $T^{big}_cM$ and, in analogy to
the case of the generalized complex structures \cite{Galt}, we
give following definition.
\begin{defin}\label{defCR} {\rm A complex, big-isotropic structure
$E\subseteq T^{big}_cM^m$ of rank $k$, with the $g$-orthogonal
subbundle $E'$, is called a {\it generalized F-structure} if
\begin{equation}\label{CRcond}E\cap\bar E'=0,\end{equation}
where the bar denotes complex conjugation. Furthermore, if $E$ is
integrable $E$ will be called a {\it generalized
CRF-structure.}}\end{defin}

If $k=m$, $E$ is a generalized (almost) complex structure. The
condition $E\cap\bar E'=0$ is equivalent with $E'\cap\bar E=0$
and, since  the orthogonal space of $\bar E$ is $\bar E'$ (because
the metric $g$ is real), $\bar E$ is a generalized F-structure too
and it is integrable iff $E$ is integrable. Furthermore, condition
(\ref{CRcond}) is equivalent to each of the equalities
\begin{equation}\label{CRdesc} T^{big}_cM=E\oplus\bar E',\;
T^{big}_cM=E'\oplus\bar E.\end{equation}

We obtain an equivalent definition of the generalized F-structures
as follows. Let $E$ be a generalized F-structure of rank $k$.
Then, $E'\cap\bar E'=S_c$ is the complexification of a real
subbundle $S\subseteq T^{big}M$ of rank
$$rank\,S=(2m-k)+(2m-k)-2m=2(m-k)$$
(use (\ref{CRdesc}) to see that $rank(E'+\bar E')=2m$) and
(\ref{CRcond}) implies $E\cap S_c=\bar E\cap S_c=0$. Thus, since
$dim(E\oplus\bar E)=2k$, decomposition (\ref{CRdesc}) becomes
\begin{equation}\label{CRdesc2} T^{big}_cM=E\oplus\bar E\oplus S_c
=L_c\oplus S_c,
\end{equation} where $L$ is the real subbundle of $T^{big}M$
such that $L_c=E\oplus\bar E$. Notice also that the definition of
$S_c$ implies $S_c\perp_g E, S_c\perp_g\bar E$.

We can use (\ref{CRdesc2}) in order to define a real endomorphism
$$\Phi:T^{big}M\rightarrow T^{big}M$$ with eigenspaces $E,\bar E,S$
of corresponding eigenvalues $\sqrt{-1},-\sqrt{-1},0$. This
endomorphism has the property
\begin{equation}\label{condYano} \Phi^3+\Phi=0.\end{equation}
Moreover, by checking on arguments in the various terms of
(\ref{CRdesc2}) while using the isotropy of $E,\bar E$ and the
relations $S_c\perp E,S_c\perp\bar E$, we see that $\Phi$ is
skew-symmetric in the sense that
\begin{equation}\label{relskew} g(\Phi(X,\alpha),(Y,\beta)) +
g((X,\alpha),\Phi(Y,\beta))=0.\end{equation}

If $\mathcal{X}=(X,\alpha)\in T^{big}M$ and
\begin{equation}\label{descX} \mathcal{X}=\mathcal{X}'+\mathcal{X}''+
\mathcal{X}''',\hspace{3mm}\mathcal{X}'\in E,\mathcal{X}''\in\bar E,
\mathcal{X}'''\in S_c,\end{equation} then
\begin{equation}\label{PhiPhi2} \Phi\mathcal{X}=\sqrt{-1}(\mathcal{X}'
-\mathcal{X}''),\;\Phi^2\mathcal{X}=-(\mathcal{X}'+\mathcal{X}''),
\end{equation} whence $ker\,\Phi=S$, $im\,\Phi=L$,
and $rank\,\Phi=2k$. Furthermore, we get the following expressions
of the projections on the terms of (\ref{CRdesc2}) and on $L$:
\begin{equation}\label{proiectorii} \begin{array}{l}
pr_{E}=-\frac{1}{2}(\Phi^2+\sqrt{-1}\Phi),\, pr_{\bar E}
=-\frac{1}{2}(\Phi^2-\sqrt{-1}\Phi),
\vspace{2mm}\\ pr_S=Id+\Phi^2,\,pr_L=-\Phi^2.
\end{array}\end{equation}
\begin{prop}\label{CRprinPhi} A generalized F-structure $E$ is
equivalent with an endomorphism $\Phi:T^{big}M
\rightarrow T^{big}M$ that satisfies the conditions
{\rm(\ref{condYano})}, {\rm(\ref{relskew})}.\end{prop}
\begin{proof} We just have deduced $\Phi$ from $E$. If we start with
$\Phi$, (\ref{condYano}) shows that the eigenvalues of $\Phi$ are
$\pm\sqrt{-1},0$, hence, $\forall x\in M$, we get a decomposition
(\ref{CRdesc2}) at $x$, where the projections are given by
(\ref{proiectorii}). In particular, the ranks of all the terms of
the derived decomposition (\ref{CRdesc2}) are lower semicontinuous
(i.e., non decreasing in a neighborhood of a point), which cannot
happen unless the ranks are constant; we will denote $k=rank\,E$.
Finally, (\ref{relskew}) implies that $E$ is big-isotropic and
$E\perp S_c, \bar E\perp S_c$, therefore, $E'=E\oplus S_c$ and
$E\cap\bar E'=0$.
\end{proof}

Condition (\ref{condYano}) is equivalent to the condition
\begin{equation}\label{condYano'} \Phi^2|_{im\,\Phi}=-Id.
\end{equation}

Furthermore, let us notice that a generalized F-structure $\Phi$
has a numerical invariant given by the negative inertia index $q$
of the metric induced by $g$ in $S$, to be called the {\it
negative index} of either $\Phi$ or $E$. Indeed, from the
definition of $S$, it follows that $g|_S$ is non degenerate (and
so is $g|_L$) and it is known that, then, $q=const.$ on the
connected components of $M$.

Accordingly, we get
\begin{prop}\label{caractcurepere} If $M$ is a connected manifold,
a skew-symmetric endomorphism $\Phi\in End(T^{big}M)$ is a
generalized F-structure of negative index $q$ iff around each
point $x\in M$ there exist local, pairwise $g$-orthogonal cross
sections
$\mathcal{Z}_a,\mathcal{Z}_\alpha$$(a=1,...,q,\alpha=1,...,p)$ of
$T^{big}M$ that satisfy the conditions
$g(\mathcal{Z}_a,\mathcal{Z}_a)=-1$,
$g(\mathcal{Z}_\alpha,\mathcal{Z}_\alpha)=+1$ and one has
\begin{equation}\label{Philocal} \Phi\mathcal{Z}_a=0,\;
\Phi\mathcal{Z}_\alpha=0,\; \Phi^2=-Id +
\sum_a(\flat_g\mathcal{Z}_\alpha)\otimes\mathcal{Z}_\alpha -
\sum_a(\flat_g\mathcal{Z}_a)\otimes\mathcal{Z}_a.\end{equation}\end{prop}
\begin{proof} The equalities
(\ref{Philocal}) imply (\ref{condYano}) and
$ker\,\Phi=span\{\mathcal{Z}_a,\mathcal{Z}_\alpha\}$, which
justifies the value of the negative index. Conversely, we know
that (\ref{condYano}) implies the existence of $S$ and the fact
that $g|_S$ is non degenerate of a negative index $q$. Then, if we
take a local basis $(\mathcal{Z}_a,\mathcal{Z}_\alpha)$ of $S$ as
required by the proposition, the conditions (\ref{Philocal}) hold.
\end{proof}

Proposition \ref{caractcurepere} and a known classical situation
\cite{{Bl},{GY}} suggest giving the following definition.
\begin{defin}\label{reperecompl} {\rm A generalized F-structure
defined by an endomorphism $\Phi\in End\,T^{big}M$ and a global,
$g$-orthogonal frame $\mathcal{Z}_a,\mathcal{Z}_\alpha$ that
satisfy all the hypotheses of Proposition \ref{caractcurepere} is
a {\it generalized F-structure with complementary
frames}.}\end{defin}

We can also derive an expression of the integrability condition of
$E$ in terms of $\Phi$. For this purpose, we consider the {\it
Courant-Nijenhuis torsion} of $\Phi$ defined by
\begin{equation}\label{torsNC} \mathcal{N}_\Phi(\mathcal{X},\mathcal{Y})
= [\Phi\mathcal{X},\Phi\mathcal{Y}]
-\Phi[\Phi\mathcal{X},\mathcal{Y}]
-\Phi[\mathcal{X},\Phi\mathcal{Y}]
+\Phi^2[\mathcal{X},\mathcal{Y}],\end{equation} where the brackets
are Courant brackets. Then we get
\begin{prop}\label{integrcuNij} The generalized F-structure $E$ is
integrable iff the corresponding endomorphism $\Phi$ satisfies one
of the following equivalent conditions
\begin{equation}\label{condCR2}
\mathcal{N}_\Phi(\mathcal{X},\mathcal{Y})=pr_S[\mathcal{X},\mathcal{Y}],
\hspace{3mm}\forall\mathcal{X},\mathcal{Y}\in\Gamma L,\end{equation}
\begin{equation}\label{condCR3}
\mathcal{S}_\Phi(\mathcal{X},\mathcal{Y})=
[\Phi\mathcal{X},\Phi\mathcal{Y}]
-[\Phi^2\mathcal{X},\Phi^2\mathcal{Y}] +
\Phi[\Phi\mathcal{X},\Phi^2\mathcal{Y}] +
\Phi[\Phi^2\mathcal{X},\Phi\mathcal{Y}]=0,\end{equation}
$\forall\mathcal{X},\mathcal{Y}\in\Gamma T^{big}M$.
\end{prop}
\begin{proof} The equivalence between the two conditions follows
by using arguments of the form $\Phi\mathcal{X},\Phi\mathcal{Y}$
in (\ref{condCR2}) (remember that $L=im\,\Phi$). Then, the
conclusion of the proposition follows from the fact that
$\forall\mathcal{X}',\mathcal{Y}'\in\Gamma
E\,\forall\mathcal{Y}''\in\Gamma\bar E$ one has
$$\mathcal{S}_\Phi(\Phi\mathcal{X}',\Phi\mathcal{Y}')=
-2([\mathcal{X}',\mathcal{Y}']+
\sqrt{-1}\Phi[\mathcal{X}',\mathcal{Y}']), \,
\mathcal{S}_\Phi(\Phi\mathcal{X}',\Phi\mathcal{Y}'')=0.$$ \end{proof}

Notice that, under the conditions (\ref{condYano}) and
(\ref{relskew}), the concomitant $\mathcal{S}_\Phi$ is
$C^\infty(M)$-linear while $ \mathcal{N}_\Phi$ is not. Another
interesting fact is given by
\begin{prop}\label{NijLS} If $E$ is an integrable big-isotropic
structure with the endomorphism $\Phi$ then
\begin{equation}\label{eqNijLS}
\mathcal{N}_\Phi(\mathcal{X},\mathcal{Y})=0,
\hspace{3mm}\forall\mathcal{X}\in L,\mathcal{Y}\in S.\end{equation}
\end{prop}
\begin{proof} The conclusion is equivalent with
\begin{equation}\label{auxNijLS}
\Phi^2[\mathcal{X},\mathcal{Y}]-\Phi[\Phi\mathcal{X},\mathcal{Y}]=0,
\hspace{3mm}\forall\mathcal{X}\in L,\mathcal{Y}\in S.
\end{equation} It suffices to check the result for $\mathcal{X}\in
E$ since the case $\mathcal{X}\in\bar E$ will be obtained by
complex conjugation. If $\mathcal{X}\in E$, (\ref{auxNijLS})
becomes
$$\sqrt{-1}\Phi[\mathcal{X},\mathcal{Y}] -
\Phi^2[\mathcal{X},\mathcal{Y}]=0,$$ which, by (\ref{proiectorii}),
is equivalent with $[\mathcal{X},\mathcal{Y}]\in E\oplus S=E'$.
For an integrable structure $E$, the previous equality holds since
the integrability of $E$ is equivalent with $[E,E']\subseteq
E'$.\end{proof}

Like in the case of the generalized almost complex structures
\cite{Galt}, one can give a representation of a generalized
CRF-structure in terms of classical tensor fields. An endomorphism
$\Phi:T^{big}M\rightarrow T^{big}M$ that satisfies the skew-symmetry
condition (\ref{relskew}) has the matrix representation
\begin{equation}\label{matriceaPhi} \Phi\left(
\begin{array}{c}X\vspace{2mm}\\ \alpha \end{array}
\right) = \left(\begin{array}{cc} A&\sharp_\pi\vspace{2mm}\\
\flat_\sigma&-\hspace{1pt}^t\hspace{-1pt}A\end{array}\right)
\left( \begin{array}{c}X\vspace{2mm}\\
\alpha \end{array}\right) \end{equation} where $(X,\alpha) \in
T^{big}M$, $A\in End(TM)$, $\sigma\in\Omega^2(M)$,
$\pi\in\chi^2(M)$, and $t$ denotes transposition. From
(\ref{matriceaPhi}) we get
\begin{equation}\label{Phi,Phi2}\begin{array}{l}
\Phi(X,\alpha)=(AX+\sharp_\pi\alpha,\flat_\sigma
X-\,^t\hspace{-1pt}A\alpha), \vspace{2mm}\\
\Phi^2(X,\alpha)=(\tilde{A}X +
\sharp_{\tilde\pi}\alpha,\flat_{\tilde\sigma}X+
\hspace{1pt}^t\hspace{-1pt}\tilde{A}\alpha),
\end{array}\end{equation}
where \begin{equation}\label{tildeA}
\tilde{A}=A^2+\sharp_\pi\flat_\sigma,\,\sharp_{\tilde\pi}=A\sharp_\pi-
\sharp_\pi\hspace{1pt}^t\hspace{-1pt}A,\,\flat_{\tilde\sigma}=
\flat_\sigma
A-\hspace{1pt}^t\hspace{-1pt}A\flat_\sigma.\end{equation} The first
formula (\ref{Phi,Phi2}) yields the following interpretation of the
entries of the matrix (\ref{matriceaPhi})
\begin{equation}\label{clasiccuPhi}
\begin{array}{l}
\pi(\alpha,\beta)=2g(\Phi(0,\alpha),(0,\beta)),\vspace{2mm}\\
\sigma(X,Y)=2g(\Phi(X,0),(Y,0)),\vspace{2mm}\\
<AX,\alpha>=2g(\Phi(X,0),(0,\alpha)).\end{array}\end{equation} The
endomorphism (\ref{matriceaPhi}) is a generalized F-structure
(i.e., satisfies (\ref{condYano})) iff:
\begin{equation}\label{condYano2}
CA=-\sharp_{\tilde\pi}\flat_\sigma,\,
C\sharp_\pi=\sharp_{\tilde\pi}\,^t\hspace{-1pt}A,\,
\flat_\sigma C=\,^t\hspace{-1pt}A\flat_{\tilde\sigma},\hspace{3mm}
C=A^2+\sharp_\pi\flat_\sigma+Id.\end{equation}

The expressions of the integrability conditions of the structure
$\Phi$ by means of the tensor fields $A,\pi,\sigma$ of
(\ref{matriceaPhi}) can be obtained by calculating (\ref{condCR3})
in each of the following cases: a) $\mathcal{X}=(0,\alpha)$,
$\mathcal{Y}=(0,\beta)$, b) $\mathcal{X}=(X,0)$,
$\mathcal{Y}=(0,\beta)$, c) $\mathcal{X}=(X,0)$,
$\mathcal{Y}=(Y,0)$. In the general case, these expressions are
complicated and give no hope for applications. For instance, the
condition
$$pr_{TM}\mathcal{S}((0,\alpha),(0,\beta))=0$$ produces the
integrability condition
\begin{equation}\label{exintegerA}
[\sharp_\pi\alpha,\sharp_\pi\beta] -
[\sharp_{\tilde\pi}\alpha,\sharp_{\tilde\pi}\beta]
+A[\sharp_\pi\alpha,\sharp_{\tilde\pi}\beta]
+A[\sharp_{\tilde\pi}\alpha,\sharp_\pi\beta]\end{equation}
$$=\sharp_\pi\{L_{\sharp_\pi\beta}(\alpha\circ\tilde{A})-
L_{\sharp_\pi\alpha}(\beta\circ\tilde{A}) +
L_{\sharp_{\tilde\pi}\beta}(\alpha\circ\tilde{A})-
L_{\sharp_{\tilde\pi}\alpha}(\beta\circ\tilde{A})$$
$$+\frac{1}{2}d(\pi(\alpha,\beta\circ\tilde{A})-\pi(\beta,\alpha\circ\tilde{A})
+\tilde{\pi}(\alpha,\beta\circ{A})-\pi(\beta,\alpha\circ{A}))\}.$$
Under the restrictive conditions $\tilde A=-Id,\tilde\pi=0$, e.g.,
in the case of a generalized complex structure \cite{Galt},
(\ref{exintegerA}) is equivalent with the fact that $\pi$ is a
Poisson bivector field.
\begin{example}\label{exgencontact} {\rm A
generalized, almost contact structure of codimension $h$
\cite{V-stable} is a system of tensor fields
$(P,\theta,F,Z_a,\xi^a)$ $(a=1,...,h)$ where $F\in End(TM)$,
$P\in\chi^2(M), \theta\in\Omega^2(M)$, $Z=(Z_a):T^*M\rightarrow
\mathds{R}^h$ is a sequence of $h$ vector fields and $\xi=(\xi^a):
TM\rightarrow\mathds{R}^h$ is a sequence of $h$ $1$-forms and the
following conditions hold
\begin{equation}\label{condF}
\begin{array}{l}P(\alpha\circ F,\beta)=P(\alpha,\beta\circ F),\;
\theta(FX,Y)=\theta(X,FY),\vspace{2mm}\\ F(Z_a)=0,\;\xi^a\circ
F=0,\;i(Z_a)\theta=0,\;i(\xi^a)P=0,\;\xi^a(Z_b)=\delta^a_b,\vspace{2mm}\\
F^2=-Id-\sharp_P\circ\flat_\theta+\sum_{a=1}^h\xi^a\otimes
Z_a.\end{array}\end{equation} Then, the tensor fields
$A=F,\pi=P,\sigma=\theta$ define a generalized F-structure $\Phi$ of
matrix (\ref{matriceaPhi}) because the conditions (\ref{condF})
imply the conditions (\ref{condYano2}). It is also easy to check
(\ref{Philocal}) with $\mathcal{Z}_\alpha=(Z_\alpha,-\xi^\alpha),
\mathcal{Z}_a=(Z_a,\xi^a)$. Hence, $\Phi$ is a generalized
F-structure with complementary frames.

The structure $(P,\theta,F,Z_a,\xi^a)$ is equivalent with the
generalized, almost complex structure defined on
$M\times\mathds{R}^h$ by the matrix
$$\Psi= \left(\begin{array}{cc} A'&\sharp_{\pi'}\vspace{2mm}\\
\flat_{\sigma'} &-\hspace{1pt}^t\hspace{-1pt}A'\end{array}\right)$$
where $$ A'=F, \;
\pi'=P+\sum_{a=1}^hZ_a\wedge\frac{\partial}{\partial t^a},\;
\sigma'=\theta+\sum_{a=1}^h\xi^a\wedge dt^a,$$ and $t^a$ are
coordinates on $\mathds{R}^h$.

Furthermore, the structure $(P,\theta,F,Z_a,\xi^a)$ is said to be
normal if $\Psi$ is integrable \cite{V-stable}. We shall prove that,
if $(P,\theta,F,Z_a,\xi^a)$ is normal, the corresponding generalized
F-structure $\Phi$ is a CRF-structure. For this purpose, we identify
$$T^{big}(M\times\mathds{R}^h)\approx
T^{big}M\oplus\mathds{R}^{2h}$$ and write $\Psi$ under the form
$$\Psi= \left(\begin{array}{cc}\Phi&\mathcal{Z}'\vspace{2mm}\\
\mathcal{Z}&0\end{array}\right),\hspace{2mm}
\mathcal{Z}=\left(\begin{array}{cc}0&Z\vspace{2mm}\\ \xi&0
\end{array}\right),\hspace{2mm} \mathcal{Z}'=
\left(\begin{array}{cc}0&-\hspace{1pt}^{t}\hspace{-1pt}Z\vspace{2mm}\\
-\hspace{1pt}^{t}\hspace{-1pt}\xi&0\end{array}\right),$$ where
$Z,\xi$ are $1$-column matrices. The integrability of $\Psi$ means
$\mathcal{N}_\Psi(\tilde{\mathcal{X}},\tilde{\mathcal{Y}})=0$,
$\forall\tilde{\mathcal{X}},\tilde{\mathcal{Y}} \in\Gamma
T^{big}(M\times\mathds{R}^h)$, where we may write
$$\tilde{\mathcal{X}}=
(X,u)\oplus(\alpha,v)=(\mathcal{X},u,v)\hspace{2mm}
(\mathcal{X}=(X,\alpha)\in\Gamma T^{big}M,\,u,v\in\mathds{R}^h)$$
and a similar expression for $\tilde{\mathcal{Y}}$. In particular,
we must have
$$\mathcal{N}_\Psi(\Psi(\mathcal{X},0,0),\Psi(\mathcal{Y},0,0))=0$$
and this equality coincides with (\ref{condCR3}). The converse may
not hold, i.e., the integrability of $\Phi$ does not imply the
normality of $(P,\theta,F,Z_a,\xi^a)$. Notice that if $\Phi$ is
defined by a normal structure $(P,\theta,F,Z_a,\xi^a)$ then $\pi$ is
a Poisson bivector field \cite{V-stable}.}\end{example}

We end this section by indicating a connection with the theory of
submanifolds of a generalized complex manifold. Let $M$ be a
submanifold of a generalized almost complex manifold $(\tilde
M,\mathcal{J})$ ($\mathcal{J}^2=-Id$,
$g(\mathcal{J}\mathcal{X},\mathcal{Y}) +
g(\mathcal{X},\mathcal{J}\mathcal{Y})=0$,
$\forall\mathcal{X},\mathcal{Y}\in\Gamma T^{big}\tilde{M}$).
\begin{defin}\label{Fsbman} {\rm $M$ is
called an {\it F-submanifold} of $\tilde{M}$ if there exists a
normal bundle $\nu M$ $(T_M\tilde M=TM\oplus\nu M)$ such that
\begin{equation}\label{eqFsbman} T^{big}M=
(T^{big}M\cap\mathcal{J}(T^{big}M))\oplus
(T^{big}M\cap\mathcal{J}(\nu^{big}M)).\end{equation}}\end{defin}

In condition (\ref{eqFsbman}) we take $T^*M=ann\,\nu M$,
$\nu^*M=ann\,TM$, $\nu^{big}M=\nu M\oplus\nu^*M$. The
$g$-skew-symmetry of $\mathcal{J}$ implies the $g$-orthogonality
of the two terms of the sum (\ref{eqFsbman}). Since the dimension
of both terms of (\ref{eqFsbman}) is upper semicontinuous (it
cannot increase in a neighborhood of $x\in M$), this dimension is
constant and the formula
\begin{equation}\label{PhiptJ} \Phi(\mathcal{X})=\left\{
\begin{array}{ll}\mathcal{J}(\mathcal{X})&\hspace{2mm}{\rm
if}\;\;\mathcal{X}\in T^{big}M\cap\mathcal{J}(T^{big}M),
\vspace{2mm}\\ 0 &\hspace{2mm}{\rm if}\;\;\mathcal{X}\in
T^{big}M\cap\mathcal{J}(\nu^{big}M)
\end{array}\right.\end{equation}
defines a generalized F-structure $\Phi$ on $M$. The structure
$\Phi$ has the associated subbundles
$L=T^{big}M\cap\mathcal{J}(T^{big}M),
S=T^{big}M\cap\mathcal{J}(\nu^{big}M)$.
\begin{defin}\label{CRFsbman} {\rm If the structure $\Phi$ defined
by formula (\ref{PhiptJ}) is CRF the submanifold $M$  will be called
a {\it CRF-submanifold} of $\tilde{M}$.}\end{defin}
\begin{rem}\label{obsBej} {\rm Definitions \ref{Fsbman},
\ref{CRFsbman} are inspired by the notion of a CR-submanifold $M$ of
a Hermitian manifold $\tilde{M}$ with the complex structure $J$
and the metric $\gamma$ \cite{Bej}, where $TM=(TM\cap J(TM))\oplus
(TM\cap J(TM)^{\perp_\gamma})$.}\end{rem}
\begin{prop}\label{condintLsbman} Any F-submanifold  $M$  of
a generalized complex manifold $(\tilde M,\mathcal{J})$ is a
CRF-submanifold.
\end{prop} \begin{proof} We shall check that the integrability condition
(\ref{condCR2}) holds. Take cross sections
$$\mathcal{X}=(X,\alpha),\mathcal{Y}=(Y,\beta),\mathcal{Z}=(Z,\gamma)
\in\Gamma T^{big}M$$ and extend $X,Y,Z,\alpha,\beta,\gamma$ to
fields $\tilde X,\tilde Y,\tilde
Z,\tilde\alpha,\tilde\beta,\tilde\gamma$ on $\tilde M$. It is easy
to check that the Courant brackets on $M$ and $\tilde M$ are
related by the following equality
\begin{equation}\label{relCMM} g_M([\mathcal{X},\mathcal{Y}]_M,\mathcal{Z})
=g_{\tilde{M}}([\tilde{\mathcal{X}},\tilde{\mathcal{Y}}]_{\tilde{M}},
\tilde{\mathcal{Z}})|_M.\end{equation} (If we take
$\mathcal{Z}\in\Gamma\nu^{big}M$ the right hand side of
(\ref{relCMM}) depends on the choice of the extensions and
(\ref{relCMM}) may not hold.)

For $\mathcal{X},\mathcal{Y}\in\Gamma(
T^{big}M\cap\mathcal{J}(T^{big}M))$ and $\mathcal{Z}\in\Gamma
T^{big}M$, using (\ref{PhiptJ}) and (\ref{relCMM}), we get
\begin{equation}\label{auxJ}
g_M(\mathcal{N}_\Phi(\mathcal{X},\mathcal{Y}),\mathcal{Z}) =
g_M(pr_S[\mathcal{X},\mathcal{Y}],\mathcal{Z})+K,\end{equation}
where
\begin{equation}\label{auxJ2}K=g_M([\Phi\mathcal{X},\Phi\mathcal{Y}]
-\Phi([\Phi\mathcal{X},\mathcal{Y}]
+[\mathcal{X},\Phi\mathcal{Y}])
-[\mathcal{X},\mathcal{Y}],\mathcal{Z})\end{equation}
$$=g_{\tilde{M}}([\mathcal{J}\tilde{\mathcal{X}},
\mathcal{J}\tilde{\mathcal{Y}}]
-\Phi([\mathcal{J}\tilde{\mathcal{X}},\tilde{\mathcal{Y}}]
+[\tilde{\mathcal{X}},\mathcal{J}\tilde{\mathcal{Y}}])
-[\tilde{\mathcal{X}},\tilde{\mathcal{Y}}],\mathcal{Z})|_M.$$
Furthermore, $\mathcal{J}$ satisfies the integrability condition
\begin{equation}\label{NJ}\mathcal{N}_{\mathcal{J}}(\mathcal{X},\mathcal{Y} =
[\mathcal{J}\tilde{\mathcal{X}}, \mathcal{J}\tilde{\mathcal{Y}}]
-\mathcal{J}([\mathcal{J}\tilde{\mathcal{X}},\mathcal{Y}]
+[\mathcal{X},\mathcal{J}\tilde{\mathcal{Y}}])
-[\mathcal{X},\mathcal{Y}]=0\end{equation} and since
$\mathcal{X},\mathcal{Y}\in\Gamma L(\Phi)\subseteq\chi^1(M)$,
(\ref{NJ}) implies
$$[\mathcal{J}\tilde{\mathcal{X}},\mathcal{Y}]
+[\mathcal{X},\mathcal{J}\tilde{\mathcal{Y}}]\in\Gamma L(\Phi),$$
therefore, the last appearance of $\Phi$ in (\ref{auxJ2}) may be
replaced by $\mathcal{J}$ and we get $K=0$. Finally, (\ref{auxJ})
with $K=0$ is equivalent with the required integrability condition.
\end{proof}
\section{CRF-structures of classical type}
In this section we discuss the simplest classes of generalized
F-structures.
\begin{defin}\label{defclasic} {\rm A generalized F-structure $\Phi$
such that $\Phi(TM)\subseteq TM$ and $\Phi(T^*M)\subseteq T^*M$ is
called a {\it classical F-structure}.}\end{defin}

The first formula (\ref{Phi,Phi2}) shows that the generalized
F-structure $\Phi$ is classical iff  $\pi=0,\sigma=0$ in the
matrix representation (\ref{matriceaPhi}). Then, (\ref{condYano2})
reduces to (\ref{condYano}) for the tensor field $F=A$ and the
generalized F-structure reduces to a Yano F-structure.

Accordingly, we may write
\begin{equation}\label{eigenF}
T_cM=H\oplus\bar H\oplus Q_c=P_c\oplus Q_c, \end{equation} where
$H,\bar H,Q_c$ are the $(\pm\sqrt{-1},0)$-eigenbundles of $F$,
respectively, $P_c=P\otimes\mathds{C}=H\oplus\bar H$, $P\subseteq
TM$ and $Q_c=Q\otimes\mathds{C}$, $Q\subseteq TM$. The terms of the
decomposition (\ref{eigenF}) have a constant dimension, because all
three dimensions are lower semicontinuous functions of $x\in M$, and
$rank\,F=const.$

It is well known that a decomposition (\ref{eigenF}) is equivalent
with the classical F-structure $F$ and the projections on the terms
of (\ref{eigenF}) are defined by the formulas (\ref{proiectorii})
with $\Phi$ replaced by $F$. Using this fact and formulas
(\ref{Phi,Phi2}), it follows that the big-isotropic complex bundle
equivalent with the generalized F-structure under discussion is
given by
\begin{equation}\label{exCR} E=H\oplus ann(H\oplus Q_c)
\end{equation}
and its orthogonal bundle is
\begin{equation}\label{exCR'} E'=(H\oplus Q_c)\oplus ann\,H.
\end{equation} Of course, we have $E\cap\bar E'=0$ and, again, we
see that a classical F-structure defines a generalized
F-structure. The subbundles $S,L$ of the structure (\ref{exCR})
are given by \begin{equation}\label{SLclasic} S=Q\oplus
ann\,P,\;L=P\oplus ann\,Q.\end{equation}
\begin{prop}\label{integrclas} The big-isotropic structure
$E$ of a classical F-structure is integrable iff
\begin{equation}\label{closureCR}
[H,H]\subseteq H,\;[H,Q_c]\subseteq H\oplus Q_c\end{equation} where
the bracket is Lie bracket on $M$.
\end{prop}
\begin{proof} The Courant bracket
(\ref{crosetC}) with $(X,\alpha),(Y,\beta)\in\Gamma E$ where $E$ is
given by (\ref{exCR}) belongs to $\Gamma E$ iff (\ref{closureCR})
holds.
\end{proof}
\begin{corol}\label{corolCRclas} For a classical F-structure, the {\it holomorphic
distribution} $H$ is a classical CR-structure.
\end{corol}
\begin{proof} A CR-structure is characterized by $H\cap\bar H=0$
plus the first condition (\ref{closureCR}) \cite{DT}.\end{proof}

Another form of the integrability condition is obtained by using the
classical Nijenhuis tensor $N_F$ given by formula (\ref{torsNC})
with arguments in $\chi^1(M)$ and Lie brackets instead of the
Courant brackets.
\begin{prop}\label{clasicCRF} The big-isotropic structure
$E$ of a classical F-structure $F$ is integrable iff
\begin{equation}\label{clasicCRF1} \begin{array}{ll}
N_F(X,Y)=pr_Q[X,Y],&\hspace{3mm}\forall X,Y\in P,\vspace{2mm}\\
N_F(X,Y)=0,&\hspace{3mm}\forall X\in P,\,Y\in
Q.\end{array}\end{equation}\end{prop}
\begin{proof} By using eigenvectors as arguments, we see that the
two conditions (\ref{clasicCRF1}) are equivalent to the two
conditions (\ref{closureCR}), respectively. \end{proof}
\begin{rem}\label{obsclasicCRF} {\rm In the case of Proposition
\ref{clasicCRF} the second condition (\ref{clasicCRF1}) is not a
consequence of the first. This is justified by the following
example. Let $H\subseteq T_cM$ be a CR-structure that is {\it
Nirenberg integrable}, i.e., $H\oplus\bar H$ is also involutive, and
let $Q$ be a complementary subbundle of $H\oplus\bar H$ in $T_cM$.
Then, it follows from the Nirenberg-Frobenius theorem \cite{Nir}
that $M$ has local complex coordinates $z^\alpha$ and real
coordinates $y^u$ such that
\begin{equation}\label{NirFrob} H=span\{\frac{\partial}{\partial
z^\alpha}\},\; Q=span\{Y_u=\frac{\partial}{\partial
y^u}+\lambda_u^\alpha\frac {\partial}{\partial z^\alpha} +
\bar{\lambda}_u^\alpha\frac {\partial}{\partial\bar
z^\alpha}\}.\end{equation} The first condition in either
(\ref{closureCR}) or (\ref{clasicCRF1}) holds but the second is
satisfied iff $\partial\lambda^\alpha_u/\partial\bar
z^\alpha=0$.}\end{rem}

We shall also notice that it is possible to express the
integrability of $E$ by a single non-skew-symmetric condition:
\begin{prop}\label{integrcuNij1}
The big-isotropic structure $E$ of a classical F-structure $F$ is
integrable iff
\begin{equation}\label{normalCR}
N_F(X,Y) = [F^2X,F^2Y] -F([X,FY] + [F^2X,FY])
\end{equation}
$$+F^2([F^2X,Y] + [F^2X,
F^2Y] + [X,Y])\hspace{2mm}(\forall X,Y\in\chi^1(M)).$$ \end{prop}
\begin{proof} Check the result by using eigenvectors as arguments.
\end{proof}

The following notion seems to be new and might be of an independent
interest.
\begin{defin}\label{defCRclas} {\rm A classical F-structure
that satisfies the conditions (\ref{closureCR}) (equivalently,
(\ref{clasicCRF1}), (\ref{normalCR})) is called a {\it classical
CRF-structure}.}\end{defin}

A classical CRF-structure may be seen either as a {\it normalized
CR-structure} $(H,Q)$ (i.e., a CR-structure $H$ with a {\it normal}
bundle $Q$ that satisfies (\ref{closureCR})) or as a {\it CR-flag}
$H\subseteq H'\subseteq T_cM$, where the following conditions are
satisfied
\begin{equation}\label{eqflag} H\cap\bar H'=0,\;H\oplus\bar
H'=T_cM,\;[H,H]\subseteq H,\;[H,H']\subseteq H'. \end{equation} A
pair $(H,Q)$ yields the flag $(H,H'=H\oplus Q)$; conversely, a
flag $(H,H')$ yields the pair $(H,Q_c=H'\cap\bar H')$.
\begin{prop}\label{clasiclocal} A tensor field $F\in End(TM)$ is an
F-structure iff around each point $x\in M$ there exist local vector
fields $Z_a$ and local $1$-forms $\xi^a$ $(a=1,...,h)$ such that
\begin{equation}\label{defF3} \xi^a(Z_b)=\delta^a_b,\;F(Z_a)=0,\;
\xi^a\circ F=0,\; F^2=-Id+\sum_{a=1}^hZ_a\otimes\xi^a.
\end{equation} Furthermore, $F$ is a CRF-structure iff the
following conditions hold
\begin{equation}\label{normal5} \mathcal{N}_F(X,Y) =
-\sum_{a=1}^hd\xi^a(X,Y)Z_a +
\sum_{a,b,c=1}^h\xi^a(X)\xi^b(Y)\xi^c([Z_a,Z_b])Z_c\end{equation}
$$+\sum_{a,b=1}^h(\xi_a(X)(L_{Z_a}\xi^b)(Y)Z_b
+\xi_a(Y)(L_{Z_a}\xi^b)(X)Z_b)$$
$$-\frac{1}{2}\sum_{a=1}^h(\xi^a(X)F(L_{Z_a}F)(Y)
-\xi^a(X)F(L_{Z_a}F)(Y)),$$
\begin{equation}\label{normal6}\sum_{a=1}^h(\xi^a(X)F(L_{Z_a}F)(Y)
+\xi^a(Y)F(L_{Z_a}F)(X))=0.\end{equation}
\end{prop}
\begin{proof} The first assertion is known and can be justified as follows.
Conditions (\ref{defF3}) imply (\ref{condYano}) for $F$. Conversely,
we know that (\ref{condYano}) implies $rank\,F=const.$ Then, if we
take an arbitrary local basis $(Z_a)_{a=1}^h$ of $Q$, the
decomposition (\ref{eigenF}) shows the existence of a unique system
of local $1$-forms $\xi^a\in ann\,L$ such that the first three
conditions (\ref{defF3}) hold. Now, we can separately check the last
condition (\ref{defF3}) on $X\in L$ and on
$Y=\sum_{a=1}^h\xi^a(X)Z_a\in Q$. For the new CRF-conditions, we use
(\ref{normalCR}), write separately the annulation of the symmetric
and skew-symmetric part, insert the local representation
(\ref{defF3}) and do the required technical computations.
\end{proof}
\begin{corol}\label{corolnormalcl}
A normal F-structure with complementary frames is a
CRF-structure.\end{corol}
\begin{proof} A normal F-structure with complementary frames
is a set of global tensor fields $(F,Z_a,\xi^a)$ $(a=1,...,h)$ that
satisfies (\ref{defF3}) and the normality condition
\begin{equation}\label{normalcontact} \mathcal{N}_F(X,Y) =
-\sum_{a=1}^hd\xi^a(X,Y)Z_a.\end{equation} It is known (e.g.,
\cite{V-stable}) that the normality condition implies the following
properties
\begin{equation}\label{auxcontact}
[Z_a,Z_b]=0,\;L_{Z_a}\xi^b=0,\;L_{Z_a}F=0. \end{equation}
Technical computations show that conditions (\ref{normalcontact}),
(\ref{auxcontact}) imply (\ref{normal5}) and
(\ref{normal6}).\end{proof}

A second class of generalized F-structures that we shall consider is
defined by
\begin{defin}\label{skewtype} {\rm A generalized F-structure $\Phi$ is said to be
{\it skew classical} if $\Phi(TM)\subseteq T^*M$ and
$\Phi(T^*M)\subseteq TM$.}\end{defin}
\begin{prop}\label{caractskewclas}
A generalized F-structures $\Phi$ is skew classical iff $A=0$ in
the corresponding matrix {\rm(\ref{matriceaPhi})}, and the
structure is fully determined either by the pair
$(\mathcal{V}=im\,\sharp_\pi,\sigma)$ where
$\sigma|_{\wedge^2\mathcal{V}}$ is non degenerate and
\begin{equation}\label{hypauxtheta} TM=\mathcal{V}\oplus
ker\,\flat_\sigma\end{equation} or by the pair
$(\Sigma=im\,\flat_\sigma,\pi)$ where $\pi|_{\wedge^2\Sigma}$ is
non-degenerate and \begin{equation}\label{hypauxtheta2}
T^*M=\Sigma\oplus ker\,\sharp_\pi.\end{equation}\end{prop}
\begin{proof} The first assertion is a straightforward consequence of the
first formula (\ref{Phi,Phi2}). Furthermore, the first formula
(\ref{Phi,Phi2}) shows that
\begin{equation}\label{eqL} L=im\,\Phi= im\,\sharp_\pi\oplus im\,\flat_\sigma,
\;S=ker\,\Phi=ker\,\flat_\sigma\oplus ker\,\sharp_\pi
\end{equation} and decomposition (\ref{CRdesc2}) implies
\begin{equation}\label{descinex2}TM=im\,\sharp_\pi\oplus ker\,\flat_\sigma,\; T^*M=
im\,\flat_\sigma\oplus ker\,\sharp_\pi.\end{equation} On the other
hand, since $\Phi^2|_L=-Id$, the second formula (\ref{Phi,Phi2})
shows that
\begin{equation}\label{auxinA=0}
(\sharp_\pi\circ\flat_\sigma)|_{im\,\sharp_\pi}=-Id,
(\flat_\sigma\circ\sharp_\pi)|_{im\,\flat_\sigma}=-Id.
\end{equation}
Therefore, $\sigma|_{\wedge^2\mathcal{V}},\pi|_{\wedge^2\Sigma}$
are non-degenerate and
\begin{equation}\label{aux2inA=0}
\sharp_\pi|_\Sigma=-(\flat_\sigma|_{\mathcal{V}})^{-1},\;
\flat_\sigma|_{\mathcal{V}}=(\sharp_\pi|_\Sigma)^{-1}.\end{equation}
The conditions (\ref{condYano2}), which reduce to
\begin{equation}\label{YanoinA=0} (Id+\flat_\sigma\sharp_\pi)\sharp_\pi=0,
\;\flat_\sigma(Id+\sharp_\pi\flat_\sigma)=0, \end{equation} are
satisfied in view of (\ref{descinex2}) and (\ref{auxinA=0}).

Now, if we have the pair $(\mathcal{V},\sigma)$ satisfying the
required conditions, the decomposition (\ref{hypauxtheta}) yields
$T^*M=\Sigma\oplus ann\,\mathcal{V}$ and $\pi$ is defined by
(\ref{aux2inA=0}) on $\Sigma$ and by $0$ on $ann\,\mathcal{V}$. A
similar procedure may be used if we start with $(\Sigma,\pi)$.
\end{proof}
\begin{prop}\label{integrskew} The skew classical, generalized
F-structure defined by a pair $(\mathcal{V},\sigma)$ that satisfies
the hypotheses of Proposition \ref{caractskewclas} is a generalized
CRF-structure iff $\mathcal{V}$ is a foliation and $\sigma$
satisfies the condition
\begin{equation}\label{integrsigma}i(X\wedge
Y)d\sigma=0, \forall X,Y\in\mathcal{V}.\end{equation}\end{prop}
\begin{proof}
From (\ref{proiectorii}) and (\ref{Phi,Phi2}) it follows that the
complex distribution of the structure $\Phi$ defined by
$(\mathcal{V},\sigma)$ is
\begin{equation}\label{EinexA0} E=\{\sharp_\pi\flat_\sigma X+
\sqrt{-1}\sharp_\pi\alpha,\flat_\sigma\sharp_\pi\alpha+\sqrt{-1}\flat_\sigma
X\}\end{equation}
$$\stackrel{(\ref{descinex2})}{=}\{(X'+\sqrt{-1}Y',
\flat_\sigma(Y'-\sqrt{-1}X'))\,/\,X',Y'\in\mathcal{V}\})\;\;
(X'=\sharp_\pi\flat_\sigma X, Y'=\sharp_\pi\alpha).$$

Formula (\ref{EinexA0}) shows that
$E=graph(\flat_{-\sqrt{-1}\sigma}|_{\mathcal{V}_c})$, which is a
situation that was discussed in Example
\ref{exgraphtheta}. The corresponding $g$-orthogonal subbundle is
\begin{equation}\label{L'theta} E'= \{(Z+\sqrt{-1}U,
\flat_\sigma(U-\sqrt{-1}Z)+(\lambda+\sqrt{-1}\mu))\end{equation}
$$ /\,Z,U\in TM,\,
\lambda,\mu\in ann\,\mathcal{V}\}$$ and, of course,
$E\cap\bar E'=0$.

Like in Example \ref{exgraphtheta}, $E$ is integrable iff
$\mathcal{V}_c$ is involutive and $$i(X\wedge
Y)d(-\sqrt{-1}\sigma)=0, \forall X,Y\in\mathcal{V}_c.$$ These
conditions are equivalent with the integrability conditions of the
structure $E_\sigma$ given by (\ref{Etheta}) with $\sigma$ instead
of $\theta$, which exactly are the conditions required by the
proposition.\end{proof}
\begin{rem}\label{obs1ptskew} {\rm
Condition (\ref{integrsigma}) is equivalent with the pair of
conditions: i) $\sigma$ induces symplectic forms on the leaves of
$\mathcal{V}$, ii) $(L_Z\sigma)|_{\mathcal{V}}=0$ for any
$\mathcal{V}$-projectable vector field $Z\in ker\,\flat_\sigma$.
Indeed, (\ref{integrsigma}) evaluated on $U\in\mathcal{V}$ is
condition i) and (\ref{integrsigma}) evaluated on $Z\in
ker\,\flat_\sigma$ is ii). Since we have the decomposition
(\ref{hypauxtheta}), we are done. It suffices to use a projectable
vector field $Z$ since (\ref{integrsigma}) is a pointwise condition
and any tangent vector at a point $x\in M$ can be extended to a
projectable vector field.}\end{rem}
\begin{example}\label{Ehresmann} {\rm
Let $\pi:M\rightarrow N$ be a symplectic fibration and let
$\mathcal{H}$ be the horizontal distribution of a symplectic
Ehresmann connection on $M$
\cite{GLS}. Then, we get a skew classical,
generalized CRF-structure $E$ associated with the pair
$(\mathcal{V},\sigma)$ where $\mathcal{V}$ is the vertical
distribution (tangent to the symplectic fibers) and $\sigma$ is
the fiber-wise symplectic form extended by $0$ on $\mathcal{H}$.
Indeed, it is known that conditions i), ii) hold in the indicated
situation \cite{GLS}. The same holds for the symplectic foliation
of a regular Poisson structure if it has a complementary
distribution whose projectable vector fields are infinitesimal
automorphisms of the Poisson structure; this may be called a
regular Poisson structure with a Poisson-Ehresmann
connection.}\end{example}

It is worth noticing that we may start with a pair
$(\mathcal{V}\subseteq TM,\theta\in\Omega^2(M))$, where
$\theta|_{\wedge^2\mathcal{V}}$ is non degenerate but
(\ref{hypauxtheta}) may not hold, and still the big-isotropic
structure $E=graph(\flat_{-\sqrt{-1}\theta}|_{\mathcal{V}_c})$ is
a generalized F-structure on $M$ (it is easy to check that
$E\cap\bar E'=0$).

The corresponding tensor fields of $E$ can be deduced as follows.
The non-degeneracy of $\theta$ on $\mathcal{V}$ implies the
existence of the decompositions
\begin{equation}\label{descpttheta}
TM=\mathcal{V}\oplus\mathcal{V}^{\perp_\theta},\;T^*M=
ann\,\mathcal{V}^{\perp_\theta}\oplus
ann\,\mathcal{V}\end{equation} and also shows that one has an
isomorphism $(\flat_\theta)|_{\mathcal{V}}:\mathcal{V}\rightarrow
ann\,
\mathcal{V}^{\perp_\theta}$. Then, if we use (\ref{descpttheta})
in the expression (\ref{perpEtheta}) of the orthogonal bundle $E'$
we see that $$E'\cap\bar E'=(\mathcal{V}^{\perp_\theta}\oplus
ann\,\mathcal{V})\otimes\mathds{C}$$ and the eigenbundle $S$ of
eigenvalue $0$ is
$$S=\mathcal{V}^{\perp_\theta}\oplus ann\,\mathcal{V}.$$

Accordingly, we get the following projection formulas
$$\begin{array}{l} pr_{E}(X,\alpha)=\frac{1}{2}(pr_{\mathcal{V}}X +
\sqrt{-1}\flat_\theta^{-1}(pr_{ann\,\mathcal{V}^{\perp_\theta}}\alpha),\,
pr_{ann\,\mathcal{V}^{\perp_\theta}}\alpha-
\sqrt{-1}\flat_\theta (pr_{\mathcal{V}}X)),\vspace{2mm}\\
pr_{\bar E}(X,\alpha)=\frac{1}{2}(pr_{\mathcal{V}}X -
\sqrt{-1}\flat_\theta^{-1}(pr_{ann\,\mathcal{V}^{\perp_\theta}}\alpha),\,
pr_{ann\,\mathcal{V}^{\perp_\theta}}\alpha+
\sqrt{-1}\flat_\theta (pr_{\mathcal{V}}X)),\vspace{2mm}\\
pr_{S_c}(X,\alpha)=(pr_{\mathcal{V}^{\perp_\theta}}X,\,
pr_{ann\,\mathcal{V}}\alpha).\end{array}$$

Now, it is easy to compute $\Phi(X,0),\Phi(0,\alpha)$, where
$\Phi$ is the equivalent endomorphism of $E$, and using
(\ref{clasiccuPhi}) we deduce
$$A=0,\;\flat_\sigma=\flat_\theta\circ pr_{\mathcal{V}},\;
\sharp_\pi=-\flat_\theta^{-1}\circ
pr_{ann\,\mathcal{V}^{\perp_\theta}}.$$ Therefore, the structure
$E$ defined by the $2$-form $\theta$ that does not satisfy
(\ref{hypauxtheta}) coincides with the structure $E$ defined by
the form
$$\sigma(X,Y)=\theta(pr_{\mathcal{V}}X,pr_{\mathcal{V}}Y)$$ for
which (\ref{hypauxtheta}) holds. The integrability conditions for
$E$ defined by $\theta$ are again those indicated in Example
\ref{exgraphtheta}.

The structures of the form
$graph(\flat_{-\sqrt{-1}\theta}|_{\mathcal{V}_c})$ considered
above extend the generalized complex structures associated with a
symplectic form \cite{Galt}.

A similar discussion applies if we start with a pair
$(\Sigma,\pi)$.  Then the generalized F-structure $\Phi$ is
equivalent with the big-isotropic structure $E_{\sqrt{-1}\pi}$
defined in Example
\ref{exgraphP} and the integrability conditions are those provided
in Example \ref{exgraphP}.\\

Now, we consider another special case:
\begin{defin}\label{defclasic2} {\rm A generalized F-structure $\Phi$
such that $\Phi^2(TM)\subseteq TM$ and $\Phi^2(T^*M)\subseteq T^*M$
is called a {\it generalized F-structure with classical
square}.}\end{defin}

The second formula (\ref{Phi,Phi2}) gives the following
characteristic properties of a structure $\Phi$ with classical
square: \begin{equation}\label{classq1}
\sharp_{\tilde{\pi}}=A\sharp_\pi-\sharp_\pi\,^t\hspace{-1pt}A=0,\;
\sharp_{\tilde{\sigma}}=\flat_\sigma
A-\,^t\hspace{-1pt}A\flat_\sigma=0\end{equation} and the conditions
(\ref{condYano2}) become
\begin{equation}\label{classq2} CA=0,\;C\sharp_\pi=0,\;\flat_\sigma
C=0,\hspace{3mm} C=A^2+\sharp_\pi\flat_\sigma+Id.\end{equation}

Furthermore, the second formula (\ref{Phi,Phi2}) implies
\begin{equation}\label{classq3}\Phi^2|_{TM}=\tilde{A},\;
\Phi^2|_{T^*M}=^t\hspace{-1pt}\tilde{A}\end{equation} and we deduce that
$$ L=im(-\Phi^2)=U\oplus U^*,$$ where
$$U=im(A^2+\sharp_\pi\flat_\sigma),\;
U^*=im(^t\hspace{-1pt}A^2+\flat_\sigma\sharp_\pi).$$ Moreover, since
the projection of $T^{big}M$ on $L$ is $-\Phi^2$, we see that
$\Pi=-A^2-\sharp_\pi\flat_\sigma$ is a projector of $TM$ onto $U$
$(\Pi^2=\Pi)$, therefore, $C$ of (\ref{classq2}) is the
complementary projector of $\Pi$ and $C^2=C$. Using $\Pi^2=\Pi$, it
is easy to check that the natural pairing between $U$ and $U^*$ is
non degenerate, hence, $U^*$ may be identified with the dual space
of $U$ and a comparison with the definition of a generalized, almost
complex structure shows that $(L,\Phi|_L)$ should be seen as a {\it
generalized, complex vector bundle} on $M$.

Conversely, if we start with an almost product structure $TM=U\oplus
V$ and a generalized complex structure $\Psi$ on the bundle
$L=U\oplus U^*$ we can define a corresponding generalized
F-structure with classical square $\Phi$ on $M$. Indeed, $\Psi$ has
a matrix (\ref{matriceaPhi}) where the entries are defined on
$U,U^*$. If these entries are extended by $0$ for any case where one
of the arguments is in $V,V^*$, the result is a matrix
(\ref{matriceaPhi}) that defines a generalized F-structure with
classical square $\Phi$. It is easy to see that the bundle $U$ of
$\Phi$ is the given one and that $\Phi|_L$ is the given structure
$\Psi$ (in particular, the projectors of the almost product
structure are $\Pi=-A^2-\sharp_\pi\flat_\sigma,
C=A^2+\sharp_\pi\flat_\sigma+Id$ because the two sides of the latter
equalities have the same values on $U$ and $V$).

We finish by recalling the notion of a $B$-field transformation:
\begin{equation}\label{Bfield} (X,\alpha)\mapsto(X,\alpha+i(X)B),
\end{equation} where $B\in\Omega^2(M)$ \cite{Galt}. Obviously, if a $B$-field
transformation is applied to a generalized F-structure, we get a
generalized F-structure again. Moreover, if $B$ is closed, the
$B$-field transformation preserves the Courant bracket, hence, if
the original structure is CRF the transformed structure is CRF as
well.
\section{Generalized metric CRF-structures}
The study of classical F-structures also includes the metric case,
i.e., the case where the structure group of $TM$ is reduced to
$U(k)\times O(h)$ \cite{{Bl},{Y}}. Equivalently, a {\it classical
metric F-structure} $(F,\gamma)$ consists of an $F$-structure $F$
and a Riemannian metric $\gamma$ on the manifold $M$ that satisfy
the following compatibility condition
\begin{equation}\label{compatFg} \gamma(X,FY)+\gamma(FX,Y)=0\hspace{3mm}
(X,Y\in\chi^1(M)).\end{equation} The pair $(F,\gamma)$ defines the
{\it fundamental $2$-form}
\begin{equation}\label{2ff} \Xi(X,Y)=\gamma(FX,Y).\end{equation}

In this section we discuss a corresponding generalized case.
Generalized Riemannian metrics were defined in \cite{Galt} as
reductions of the structure group $O(m,m)$ of $(T^{big}M,g)$ to
$O(m)\times O(m)$. For the reader's convenience, we recall the basic
facts given in \cite{Galt}.

A {\it generalized, Riemannian metric} is a Euclidean (positive
definite) metric $G$ on the bundle $T^{big}M$, which is compatible
with the metric $g$ given by (\ref{gFinC}) in the sense that the
musical isomorphism
\begin{equation}\label{sharpbig} \sharp_G:T^{big}M=TM\oplus T^*M \rightarrow
T^*M\oplus TM\approx T^{big}M,\end{equation} where $\approx$ is
the isomorphism $(\alpha,X)\leftrightarrow(X,\alpha)$, satisfies
the conditions
\begin{equation}\label{Riemannbig1}
\sharp_G^2=Id,\end{equation}
\begin{equation}\label{Riemannbig2}g(\sharp_G(X,\alpha),\sharp_G(Y,\beta)) =
g((X,\alpha),(Y,\beta)).
\end{equation}

The isomorphism $\sharp_G$ is determined by the formula
\begin{equation}\label{eqGg}
2g(\sharp_G(X,\alpha),(Y,\beta)) = G((X,\alpha),(Y,\beta)),
\end{equation} and, if we ask (\ref{Riemannbig1}) to hold, we see that
$\sharp_G$ may be represented in the matrix form
\begin{equation}\label{matriceaG} \sharp_G\left(\begin{array}{c}
X\vspace{2mm}\\ \alpha\end{array}\right)=
\left(\begin{array}{cc}\varphi&\sharp_\gamma\vspace{2mm}\\
\flat_\beta&^t\hspace{-1pt}\varphi\end{array}\right)
\left(\begin{array}{c} X\vspace{2mm}\\ \alpha\end{array}\right),
\end{equation} where $\varphi\in End(TM)$ and
$\beta,\gamma$ are classical Riemannian metrics on $M$.

Furthermore, condition (\ref{Riemannbig2}) is equivalent with
\begin{equation}\label{GcuRbig1} \begin{array}{c}
\varphi^2=Id-\sharp_\gamma\circ\flat_\beta,\;\gamma(\varphi X,Y)+
\gamma(X,\varphi Y)=0,\vspace{2mm}\\ \beta(\varphi
X,Y)+\beta(X,\varphi Y)=0.\end{array}\end{equation} Since
$\beta,\gamma$ are non degenerate, the first condition
(\ref{GcuRbig1}) yields
\begin{equation}\label{beta-gamma}
\flat_\beta=\flat_{\gamma}\circ(Id-\varphi^2), \;
\sharp_\gamma=(Id-\varphi^2)\circ\sharp_\beta\end{equation} and we see
that the generalized Riemannian metrics of $M$ are in a bijective
correspondence with pairs $(\gamma,\varphi)$ or $(\beta,\varphi)$
where $\gamma,\beta$ are Riemannian metrics and $\varphi$ is a
$\gamma,\beta$-skew-symmetric $(1,1)$-tensor field.

Notice that, modulo (\ref{GcuRbig1}), the fact that $\gamma$ is
positive definite implies that $\beta$ and $G$ are positive
definite; indeed, (\ref{GcuRbig1}), (\ref{beta-gamma}) and
(\ref{eqGg}) imply
$$\beta(X,X)=\gamma(X,X)+\gamma(\varphi X,\varphi X),$$ $$G((X,\flat_\gamma Y),
(X,\flat_\gamma Y))= \gamma(X,X)+\gamma(Y+\varphi X,Y+\varphi X).$$
Conversely, if $\beta$ is positive definite so is $\gamma$ too. In
particular, if $\varphi=0$, $G$ reduces to the classical Riemannian
metric $\gamma$.

Furthermore, $\varphi$ may be replaced by the $2$-form $\psi$
defined by
\begin{equation}\label{psi} \flat_\psi=-\flat_\gamma\circ\varphi,
\end{equation}
which means that one has a bijective correspondences $G
\leftrightarrow(\gamma,\psi)$.

From (\ref{Riemannbig1}), (\ref{Riemannbig2}) and (\ref{eqGg}), it
follows that a generalized, Riemannian metric $G$ produces a
decomposition
\begin{equation}\label{descompG} T^{big}M=V_+\oplus
V_-,\end{equation} where $V_\pm$ are the $(\pm1)$-eigenbundles,
which simultaneously are $G$ and $g$ orthogonal. On $V_\pm$ one has
$G=\pm 2g$, respectively, hence, $g$ is positive definite on $V_+$
and negative definite on $V_-$, whence, $dim\,V_\pm= m$. Conversely,
a decomposition (\ref{descompG}) with $m$-dimensional,
$g$-orthogonal, terms that are $g$-positive and $g$-negative,
respectively, defines the generalized, Riemannian metric
$G=2g_{V_+}-2g_{V_-}$. This exactly means that the structure group
of $T^{big}M$ was reduced to $O(m)\times O(m)$.

The projectors associated with the decomposition (\ref{descompG})
are given by
\begin{equation}\label{prpeE}
pr_{\pm}=\frac{1}{2}(Id\pm\sharp_G),\end{equation} and, if we apply
them to pairs $(0,\alpha),(X,0)$ using (\ref{matriceaG}), we see
that the projections
\begin{equation}\label{tau}
\tau_{\pm}=pr_{TM}|_{V_\pm},\tau^*_{\pm}=pr_{T^*M}|_{V_\pm}\end{equation}
are surjective hence, isomorphisms. From (\ref{matriceaG}) and
(\ref{psi}), we get
\begin{equation}\label{formuletau} \begin{array}{l}\tau_+^{-1}(X)=
(X,\flat_\gamma(X-\varphi
X))=(X,\flat_{\psi+\gamma}X),\vspace{2mm}\\
\tau_-^{-1}(X)= (X,-\flat_\gamma(X+\varphi
X))=(X,\flat_{\psi-\gamma}X),\end{array}\end{equation} therefore,
\begin{equation}\label{exprEpm}
V_\pm=\{(X,\flat_{\psi\pm\gamma}X)\,/\,X\in TM\}.\end{equation}

The isomorphisms $\tau_{\pm}$ may be used in order to transfer the
metrics $G|_{V_\pm}$ to metrics $G_\pm$ of the tangent bundle $TM$
given by
\begin{equation}\label{metriciinduse}
G_{\pm}(X,Y)=G(\tau_{\pm}^{-1}(X),\tau_{\pm}^{-1}(Y))= \pm
2g(\tau^{-1}_{\pm}(X),\tau^{-1}_{\pm}(Y))=2\gamma(X,Y),
\end{equation}
where $\gamma$ is the metric that appears in
(\ref{matriceaG}).\\

The following definition extends the one given in \cite{Galt} for
the generalized complex case.
\begin{defin}\label{defPhimetric} {\rm A {\it generalized metric
F-structure} is a pair $(\Phi,G)$, where $\Phi$ is a generalized
F-structure and $G$ is a generalized Riemannian metric, such that
the following skew-symmetry condition holds
\begin{equation}\label{PhiG} G(\Phi\mathcal{X},\mathcal{Y})
+G(\mathcal{X},\Phi\mathcal{Y})=0\hspace{3mm}(
\mathcal{X},\mathcal{Y}\in\Gamma T^{big}M).\end{equation}
}\end{defin}

Using the $g$-skew-symmetry (\ref{relskew}) of $\Phi$ and formula
(\ref{eqGg}), we see that (\ref{PhiG}) is equivalent with the
commutation condition
\begin{equation}\label{compatGPhi}
\sharp_G\circ\Phi=\Phi\circ\sharp_G. \end{equation} Condition
(\ref{compatGPhi}) implies that the pair
\begin{equation}\label{adouastr}(\Phi^c=\sharp_G\circ\Phi,G)\end{equation}
is a second generalized metric F-structure that commutes with
$\Phi$. We will refer to $\Phi^c$ as the {\it complementary
structure}. In the generalized almost complex case a commuting pair
$(\Phi,\Phi^c)$ defines $G$ by $\sharp_G=-\Phi\circ\Phi^c$
\cite{Galt}.

Let us assume that the structure $\Phi$ has the matrix
representation (\ref{matriceaPhi}). The following proposition
expresses the compatibility between $\Phi$ and $G$ via the
matrices (\ref{matriceaG}) and (\ref{matriceaPhi}).
\begin{prop}\label{metricaHermite} The pair $(G,\Phi)$, where $G$
is a generalized, Riemannian metric given by {\rm(\ref{matriceaG})}
and $\Phi$ is a generalized F-structure given by
{\rm(\ref{matriceaPhi})}, is a generalized metric F-structure iff
the following two conditions hold
\begin{equation}\label{H1} \gamma(AX,Y) + \gamma(X,AY) =
\varpi(\varphi X,Y)-\varpi(X,\varphi Y),\end{equation}
\begin{equation}\label{H2} \sigma(X,Y)=\varpi(X,Y)
- \varpi(\varphi^2X,Y)+\gamma([A,\varphi](X),Y), \end{equation}
where $[A,\varphi]=A\circ\varphi-\varphi\circ A$ and
$\varpi=\flat_\gamma\pi$ is defined by
\begin{equation}\label{defvarpi} \varpi(X,Y)=\pi(\flat_\gamma
X,\flat_\gamma Y).\end{equation} \end{prop}
\begin{proof} The commutation condition
(\ref{compatGPhi}) is equivalent with
\begin{equation}\label{compat1} \begin{array}{l} \varphi\circ
A+\sharp_\gamma\circ\flat_\sigma= A\circ
\varphi+\sharp_\pi\circ\flat_\beta,
\vspace{2mm}\\
\varphi\circ\sharp_\pi -
\sharp_\gamma\circ\,^t\hspace{-1pt}A =
A\circ\sharp_\gamma+\sharp_\pi\circ\,^t\hspace{-1pt}\varphi,\vspace{2mm}\\
\sharp_\beta\circ A+
^t\hspace{-1pt}\varphi\circ\flat_\sigma = \flat_\sigma\circ\varphi
- ^t\hspace{-1pt}A\circ\flat_\beta,\vspace{2mm}\\
\flat_\beta\circ\sharp_\pi -
^t\hspace{-1pt}\varphi\circ ^t\hspace{-1pt}A =
\flat_\sigma\circ\sharp_\gamma -
^t\hspace{-1pt}A\circ\,^t\hspace{-1pt}\varphi.\end{array}\end{equation}

Furthermore, the last condition (\ref{compat1}) is the
transposition of the first, and the first condition implies the
equivalence between the second and third condition. Indeed, the
first condition (\ref{compat1}) is equivalent with
\begin{equation}\label{determinsigma} \flat_\sigma=
\flat_\gamma\circ(A\circ\varphi-\varphi\circ A+\sharp_\pi\circ\flat_\beta).
\end{equation} If this expression of $\flat_\sigma$ is inserted in
the third condition (\ref{compat1}), while taking into account the
$\beta,\gamma$-skew-symmetry of $\varphi$ and (\ref{beta-gamma}),
and the result is composed by $\sharp_\gamma$ at the left and by
$\sharp_\beta$ at the right, one gets
\begin{equation}\label{aux3-2}\varphi\circ\sharp_\pi-
\sharp_\gamma\circ\,^t\hspace{-1pt}A =
A\circ(Id-\varphi^2)\circ\sharp_\beta +
\sharp_\pi\circ\,^t\hspace{-1pt}\varphi.\end{equation} Then, the first
condition (\ref{GcuRbig1}) shows that (\ref{aux3-2}) coincides
with the second condition (\ref{compat1}).

Thus, the compatibility between $G$ and $\Phi$ is equivalent with
(\ref{determinsigma}) together with the second condition
(\ref{compat1}). If (\ref{determinsigma}) is evaluated on tangent
vectors $X,Y$ and the second condition (\ref{compat1}) is
evaluated on $\flat_\gamma X,\flat_\gamma Y$ the required
conclusion is obtained.\end{proof}

We proceed by a recall of Gualtieri's expression of a generalized,
metric, almost complex structure by classical structures while
replacing the structure by a generalized, metric F-structure
\cite{Galt}.

From (\ref{descompG}), it follows that (\ref{compatGPhi}) is
equivalent with the $\Phi$-invariance of the subbundles $V_\pm$.
Thus, if we also take into account the skew-symmetry
(\ref{relskew}) of $\Phi$, a $G$-compatible, generalized metric
F-structure is equivalent with a pair of bundle morphisms
$\Phi_{\pm}\in End\,V_\pm$ which satisfy the condition
$\Phi_\pm^3+\Phi_\pm=0$ and are skew-symmetric with respect to
$G|_{V_{\pm}}$. Furthermore, we have decompositions
\begin{equation}\label{descVpm} V_\pm=E_\pm\oplus\bar E_\pm\oplus S_\pm
\end{equation} where
the terms are the $(\pm\sqrt{-1},0)$-eigenbundles of $\Phi_\pm$.
Hence, the eigenbundles of $\Phi$ are
\begin{equation}\label{desceigen} E=E_+\oplus E_-,\, \bar E=\bar E_+
\oplus\bar E_-,\, S=S_+\oplus S_-\end{equation}
and we have
\begin{equation}\label{descVpm2} E_\pm=V_\pm\cap E,\,
\bar E_\pm=V_\pm\cap\bar E,\,S_\pm=V_\pm\cap S.\end{equation}

Similar decompositions hold for the complementary structure
$\Phi^c$. We shall denote the corresponding vector bundles by
means of an upper index $c$ and, if we look at the corresponding
eigenvalue and use (\ref{adouastr}), we get
\begin{equation}\label{EEc}
E_+\subseteq E^c,\,E_-\subseteq\bar E^c,\,S^c=S.\end{equation}
Formulas (\ref{descVpm}), (\ref{desceigen}), (\ref{descVpm2}) hold
for $E,\bar E$ replaced by $E^c,\bar E^c$. Moreover, using
(\ref{EEc}) it is easy to get \begin{equation}\label{EEc2}
E_+=E\cap E^c,\,E_-=E\cap\bar E^c.\end{equation} Finally, notice
that (\ref{adouastr}) implies $(\Phi^c)^c=\Phi$, therefore if we
consider the decomposition (\ref{descVpm}) for $\Phi_c$ instead of
$\Phi$ we get
\begin{equation}\label{c} E_+^c=E_+,\,E_-^c=\bar
E_-,\,E^c=E_+\oplus\bar E_-.\end{equation}

Furthermore, the structures $\Phi_{\pm}$ may be transferred to
$TM$ by
\begin{equation}\label{strF}
F_{\pm}=\tau_{\pm}\circ\Phi_\pm\circ\tau_{\pm}^{-1}\in End\,TM
\end{equation} and, also recalling formula (\ref{metriciinduse}),
the conclusion is that the $G$-compatible, generalized F-structure
$\Phi$ is equivalent with the pair of classical F-structures
$F_{\pm}$ of $TM$ that satisfy the skew-symmetry condition
\begin{equation}\label{compatcutheta} \gamma(F_{\pm}X,Y) +\gamma
(X,F_{\pm}Y)=0.\end{equation}

It is easy to find the connection between $F_{\pm}$ and the matrix
(\ref{matriceaPhi}) of $\Phi$. Using (\ref{formuletau}), it
follows that
\begin{equation}\label{formuleJ}
F_{\pm}=A+\sharp_\pi\circ\flat_{\psi\pm\gamma}.
\end{equation} Conversely, from (\ref{formuleJ}) we get
\begin{equation}\label{eqpiAdinF}
\begin{array}{c}
\sharp_\pi=\frac{1}{2}(F_+-F_-)\circ\sharp_\gamma,\vspace{2mm}\\
A=\frac{1}{2}[F_+\circ(Id-\sharp_\gamma\flat_\psi)+
F_-\circ(Id+\sharp_\gamma\flat_\psi)].\end{array}
\end{equation} The remaining entry of the matrix
(\ref{matriceaPhi}) is $\flat_\sigma$ is given by
(\ref{determinsigma}).

Therefore, we have the same result as in \cite{Galt}:
\begin{prop}\label{prGaltK} A generalized metric F-structure
$(G,\Phi)$ is equivalent with a quadruple $(\gamma,\psi,F_+,F_-)$,
where $\gamma$ is a classical, Riemannian metric on $M$, $\psi$ is
a $2$-form, and $(F_{\pm},\gamma)$ are classical metric
F-structures of $M$.
\end{prop} \begin{rem}\label{FPhic} {\rm Assume
that the generalized metric F-structure $(G,\Phi)$ has the
corresponding quadruple $(\gamma,\psi,F_+,F_-)$. Since the
complementary structure $\Phi^c$ satisfies the conditions
$\Phi^c_{_\pm}=\pm\Phi_{_\pm}$, formula (\ref{strF}) shows that
$(G,\Phi^c)$ has the corresponding quadruple
$(\gamma,\psi,F_+,-F_-)$.}\end{rem}
\begin{example}\label{Fmetricclasic} {\rm Let $(F,\gamma)$ be a
classical metric F-structure, $\Phi$ the corresponding generalized
F-structure given by ({\ref{exCR}}) and $G$ the generalized
Riemannian metric defined by a classical Riemannian metric
$\gamma$ (with $\varphi=0$). Then, Proposition
\ref{metricaHermite} shows that $(\Phi,G)$ is a generalized metric
F-structure. By (\ref{psi}), this structure has the $2$-form
$\psi=0$, by (\ref{exprEpm}), $V_\pm=\{(X,\flat_{\pm\gamma}X)\}$,
by (\ref{metriciinduse}) the metric to be considered on $TM$ is
$2\gamma$ and by (\ref{strF}) we have $F_+=F_-=F$. Furthermore,
with the notation of the first part of Section 3, we have
\begin{equation}\label{eqinexFclas}
S_\pm=\{(X,\flat_{\pm\gamma})\,/\,X\in Q\},\;
E_\pm=\{(X,\flat_{\pm\gamma})\,/\,X\in
H\}.\end{equation}}\end{example}

Now, if $\Phi$ is  a generalized CRF-structure that is
skew-symmetric with respect to the generalized Riemannian metric
$G$ then $(\Phi,G)$ is a {\it generalized metric CRF-structure}.
We will extend the notion of a generalized K\"ahler manifold
\cite{Galt} by means of the following definition.
\begin{defin}\label{defCRFK} {\rm A generalized, metric F-structure
$(\Phi,G)$ with the associated eigenbundles $(E_\pm,S_\pm)$ is a
{\it generalized CRFK-structure} (and $(M,\Phi,G)$ is a {\it
generalized CRFK manifold}) iff the following Courant bracket
closure conditions hold:
\begin{equation}\label{CRFK1'} \begin{array}{l}[E_+,E_+]
\subseteq E_+,\,[E_+,S_+]\subseteq E_+ \oplus S_+,\vspace{2mm}\\

[E_-,E_-]\subseteq E_-,\, [E_-,S_-]\subseteq E_- \oplus S_-.
\end{array}\end{equation}}\end{defin}
The label K used above comes from the name of K\"ahler.

The relation between this definition and the definition of a
generalized K\"ahler manifold given in \cite{Galt} is shown by the
following proposition.
\begin{prop}\label{CRFK0} The generalized, metric F-structure
$(\Phi,G)$ is a CRFK-structure iff $\Phi$ and its complementary
structure $\Phi^c$ are CRF-structures and
\begin{equation}\label{SpSm} [S_+,S_-]\subseteq S.\end{equation}
\end{prop}
\begin{proof} We shall
use the following property of the Courant bracket (\cite{LWX},
axiom (v) of Courant algebroids):
\begin{equation}\label{axv}\begin{array}{ll}
X(g(\mathcal{Y},\mathcal{Z}))&=g([\mathcal{X},
\mathcal{Y}],\mathcal{Z}) + g(\mathcal{Y},[\mathcal{X},\mathcal{Z}])
\vspace{2mm}\\
&+\frac{1}{2}(Z(g(\mathcal{X},\mathcal{Y}) +
Y(g(\mathcal{X},\mathcal{Z}))),\end{array}\end{equation}
$\forall\mathcal{X}=(X,\alpha),\,\mathcal{Y}=(Y,\beta),\,\mathcal{Z}=(Z,\gamma)
\in\Gamma T^{big}M$.

The integrability (CRF condition) of $\Phi,\Phi^c$ means
$[E,E]\subseteq E,\,[E^c,E^c]\subseteq E^c$. From
(\ref{desceigen}), (\ref{EEc2}) and (\ref{c}), it follows
straightforwardly that these integrability conditions are
equivalent with the conditions:
\begin{equation}\label{CRFK1} [E_+,E_+]
\subseteq E_+,\, [E_-,E_-]\subseteq E_-,\, [E_+,E_-] \subseteq E,
[E_+,\bar E_-] \subseteq E^c.
\end{equation}

We will show that conditions (\ref{CRFK1}) are equivalent with the
conditions
\begin{equation}\label{CRFKp} [E_+,E_+]
\subseteq E_+,\, [E_-,E_-]\subseteq E_-,\, [E_+,E_-] \perp_g S,
[E_+,\bar E_-] \perp_g S.
\end{equation} Indeed, the third and fourth condition
(\ref{CRFK1}) imply the third and fourth condition (\ref{CRFKp}),
respectively, since $E\perp_g S,E^c\perp_g S^c=S$. Conversely, if
the first two conditions of (\ref{CRFK1}) hold and if we take $
\mathcal{X},\mathcal{Z}\in E_\pm,\mathcal{Y}\in E_\mp$ in
(\ref{axv}), we get
\begin{equation}\label{pmperpg} [E_+,E_-]\subseteq
E^{\perp_g}=E'=E\oplus S.\end{equation} Since, by the third
condition (\ref{CRFKp}), $[E_+,E_-]\subseteq E\oplus\bar E$, it
follows that
$$[E_+,E_-]\subseteq(E\oplus S)\cap(E\oplus\bar E)=E.$$ With the
same argument for $E^c$ instead of $E$ we see that the last
condition (\ref{CRFKp}) implies the last condition (\ref{CRFK1}).

Furthermore, we show that the conditions (\ref{CRFKp}) are
equivalent with either
\begin{equation}\label{CRFK2}
[E_+,E_+]\subseteq E_+,\,[E_-,E_-]\subseteq E_-,\,[E_+,S]\subseteq
E_+ \oplus S\end{equation} or
\begin{equation}\label{CRFK3}
[E_+,E_+]\subseteq E_+,\,[E_-,E_-]\subseteq E_-,\,[E_-,S]\subseteq
E_- \oplus S.\end{equation} Indeed, let us use (\ref{axv}) with
$\mathcal{X}\in E_\pm,\mathcal{Y}\in E_\mp,\mathcal{Z}\in S$. The
result is that
\begin{equation}\label{auxindem'} [E_+,E_-]\perp_g
S\,\Leftrightarrow\, [E_+,S]\perp_g E_-\,\Leftrightarrow
[E_-,S]\perp_g E_+.\end{equation} Similarly, but changing $E_-$ to
$\bar E_-$, we get
\begin{equation}\label{auxindem'-} [E_+,\bar E_-]\perp_g
S\,\Leftrightarrow\,[\bar E_-,S]\perp_g E_+\,\Leftrightarrow
\,[E_+,S]\perp_g\bar E_-.\end{equation}

The relations (\ref{auxindem'}) and (\ref{auxindem'-}) show that
(\ref{CRFKp}) implies
$$[E_+,S]\perp_g(E_-\oplus\bar E_-),
\,[E_-,S]\perp_g(E_+\oplus\bar E_+).$$ Accordingly, and since
(\ref{CRFKp}) is equivalent to the integrability of $E$, we get
$$\begin{array}{c}[E_+,S]\subseteq(E\oplus S)\cap
(E_+\oplus\bar E_+\oplus S)=E_+\oplus S,\vspace{2mm}\\

[E_-,S]\subseteq(E\oplus S)\cap(E_-\oplus\bar E_-\oplus
S)=E_-\oplus S,
\end{array} $$

\noindent which shows that (\ref{CRFKp}) implies both (\ref{CRFK2}) and
(\ref{CRFK3}). Conversely, the condition $[E_+,S]\subseteq E_+
\oplus S$ implies both $[E_+,S]\perp_g E_-$ and
$[E_+,S]\perp_g\bar E_-$ and (\ref{auxindem'}), (\ref{auxindem'-})
show that the conditions (\ref{CRFK2}) imply (\ref{CRFKp}).
Similarly, $[E_-,S]\subseteq E_-
\oplus S$ implies both $[E_-,S]\perp_g E_+$ and
$[E_-,S]\perp_g\bar E_+$, equivalently, $[\bar E_-,S]\perp_g E_+$,
and (\ref{auxindem'}), (\ref{auxindem'-}) show that the conditions
(\ref{CRFK3}) imply (\ref{CRFKp}).

Finally, if we use (\ref{axv}) with $\mathcal{X}\in
S_\pm,\mathcal{Y}\in E_\mp,\mathcal{Z}\in S_\mp$ we get
\begin{equation}\label{auxfinal}
[S_+,S_-]\perp_gE_+\,\Leftrightarrow\,[E_+,S_+]\perp_g S_-,\,
[S_+,S_-]\perp_gE_-\,\Leftrightarrow\,[E_-,S_-]\perp_g S_+.
\end{equation} Therefore, the addition of hypothesis (\ref{SpSm})
to (\ref{CRFK2}), (\ref{CRFK3}) leads to
$$\begin{array}{l}[E_+,S_+]\subseteq(E_+\oplus
S)\cap(V_+\oplus E_-\oplus\bar E_-)=E_+\oplus S_+,\vspace{2mm}\\

[E_-,S_-]\subseteq(E_-\oplus S)\cap(V_-\oplus E_+\oplus\bar
E_+)=E_-\oplus S_-.\end{array}$$

\noindent Conversely, from (\ref{auxfinal}) we get
$$\begin{array}{l}[E_+,S_+]\subseteq E_+\oplus S_+\,\Rightarrow\,
[E_+,S_+]\perp_gS_-\,\Rightarrow\,[S_+,S_-]\perp_gE_+,\vspace{2mm}\\

[E_-,S_-]\subseteq E_-\oplus S_-\,\Rightarrow\, [E_-,S_-]\perp_g
S_+\,\Rightarrow\,[S_+,S_-]\perp_gE_-.\end{array}$$

\noindent Since $S_\pm$ are real subbundles, by complex
conjugation, we also get $[S_+,S_-]\perp_g\bar
E_+,[S_+,S_-]\perp_g\bar E_-$ and (\ref{SpSm}) follows.
\end{proof}
\begin{rem}\label{comparcuGalt} {\rm During the proof of Proposition
\ref{CRFK0} we have obtained several characterizations of the integrability of the
pair of structures $(\Phi,\Phi^c)$: (\ref{CRFK1}), (\ref{CRFKp}),
(\ref{CRFK2}), (\ref{CRFK3}). In the generalized K\"ahler case
$S=0$, condition (\ref{SpSm}) does not appear and (\ref{CRFKp})
shows that the last two conditions (\ref{CRFK1}) are superfluous.
This is a simple proof of results that appeared in Proposition
6.10 and Theorem 6.28 of \cite{Galt}.}\end{rem}

Like in the case of the generalized K\"ahler structures
\cite{Galt}, it is possible to obtain conditions that are
equivalent with (\ref{CRFK1'}) and are expressed in terms of the
projected structures $F_\pm$ given by (\ref{strF}).
\begin{prop}\label{propCRFK} The generalized metric CRF-structure
$(\Phi,G)$ with the associated structures $(F_\pm,\gamma,\psi)$ is
a CRFK-structure iff $F_\pm$ are classical metric CRF-structures
and the equalities
\begin{equation}\label{CRFK4} i(X\wedge Y)d\psi= \pm (i(X)L_Y\gamma
- L_Xi(Y)\gamma)
\end{equation}  hold for either $X,Y\in H_\pm$ or $X\in H_\pm,Y\in
Q_\pm$.\end{prop}
\begin{proof}
The restriction of the Courant bracket to the subbundles $V_\pm$
defined by (\ref{formuletau}) is given by the formula (see also
\cite{Galt}):
\begin{equation}\label{auxgeneral} [(X,\flat_{\psi\pm\gamma}X),
(Y,\flat_{\psi\pm\gamma}Y)] = ([X,Y],\flat_{\psi\pm\gamma}[X,Y]
\end{equation} $$+i(X\wedge Y)d\psi
\pm (L_Xi(Y)\gamma-i(X)L_Y\gamma))\hspace{2mm}(X,Y\in\chi^1(M)).$$ This
formula follows from the general expression (\ref{crosetC}) of the
Courant bracket by evaluating the $1$-form component on a vector
field $Z$.

Let us denote by $H_\pm,\bar H_\pm,Q_\pm$ the
$(\pm\sqrt{-1},0)$-eigenbundles of $F_\pm$. From
(\ref{auxgeneral}), it follows that the CRFK-conditions
(\ref{CRFK1'}) are equivalent with the conditions
\begin{equation}\label{metricjos}
[H_\pm,H_\pm]\subseteq H_\pm,\, [H_\pm,Q_\pm]\subseteq H_\pm\oplus
Q_\pm\end{equation} together with the equalities (\ref{CRFK4}).
\end{proof}

In what follows we produce some more equivalent CRFK-conditions.
Any connection on the principal bundle of frames of $TM$ given by
the reduction of the structure group to $U(rank\,H_\pm)\times
O(rank\,Q_\pm)$ defined by the metric F-structures
$(F_\pm,\gamma)$, will be called an {\it adapted connection}. (Of
course, we have different plus-adapted and minus-adapted
connections and all the formulas where we write a double sign
$\pm$ include two different formulas.) The parallel translations
of adapted connections preserve the structure $(F_\pm,\gamma)$ and
the associated covariant derivative $\nabla^\pm$ is characterized
by the conditions
\begin{equation}\label{nablapm} \nabla^\pm\gamma=0,\; \nabla^\pm
F_\pm=0.\end{equation} Furthermore, the {\it difference tensor}
\begin{equation}\label{diftensor}
\Theta^\pm(X,Y)=\nabla^\pm_XY-\nabla_XY,\end{equation} where $\nabla$ is the
Levi-Civita connection of the metric $\gamma$ will be called the
{\it Levi-Civita difference} of the adapted connection. Since
$\nabla\gamma=0$, we must have
\begin{equation}\label{condTheta} \gamma(\Theta^\pm(X,Y),Z)+
\gamma(Y,\Theta^\pm(X,Z))=0.\end{equation} On the other hand, the condition
$\nabla^\pm F_\pm=0$ is equivalent with
\begin{equation}\label{condTheta2}
\Theta^\pm(X,F_\pm Y)-F_\pm\Theta^\pm(X,Y)=-(\nabla_XF_\pm)(Y),\end{equation}
which also implies
\begin{equation}\label{condTheta3}
\Theta^\pm(X,F_\pm^2 Y)-F_\pm^2\Theta^\pm(X,Y)=-(\nabla_XF^2_\pm)(Y).\end{equation}
Thus, the adapted connections are obtained from the Levi-Civita
connection by the addition of a difference tensor that satisfies
conditions (\ref{condTheta}) and (\ref{condTheta2}).
\begin{prop}\label{caractcuconex} Let $(\Phi,G)$ be a
generalized metric CRF-structure with the associated structures
$(F_\pm,\gamma,\psi)$ and let $\nabla^\pm$ be adapted connections
with the Levi-Civita differences $\Theta^\pm$. Then, $(\Phi,G)$ is
a CRFK-structure iff $F_\pm$ are classical metric CRF-structures
and the equalities
\begin{equation}\label{equivCRFK2}
\gamma(\Theta^\pm(Z,Y),X)=\mp\frac{1}{2}
d\psi(X,Y,Z)\end{equation} hold for any $Z\in\chi^1(M)$ and either
$X,Y\in H_\pm$ or $X\in H_\pm,Y\in Q_\pm$.\end{prop}
\begin{proof}
By a simple computation and using $\nabla^\pm\gamma=0$ we get
\begin{equation}\label{Liecutors} \begin{array}{r}(L_X\gamma)(Y,Z)=
\gamma(\nabla^\pm_YX,Z)
+\gamma(Y,\nabla^\pm_ZX)\vspace{2mm}\\ +\gamma(T^\pm(X,Y),Z)
+\gamma(Y,T^\pm(X,Z)),\end{array}\end{equation} where $T^\pm$ is
the torsion of $\nabla^\pm$. Then, if we evaluate (\ref{CRFK4}) on
$Z\in\chi^1(M)$ and use (\ref{nablapm}) and (\ref{Liecutors}) we
get the following equivalent form of (\ref{CRFK4}):
\begin{equation}\label{equivCRFK} \begin{array}{l}
d\psi(X,Y,Z)=\pm [\gamma(X,\nabla^\pm_ZY)-\gamma(Y,\nabla^\pm_ZX)
\vspace{2mm}\\ +\gamma(X,T^\pm(Y,Z)) +\gamma(Y,T^\pm(Z,X))
-\gamma(Z,T^\pm(X,Y)),\end{array}\end{equation} where the first
two terms of the right hand side vanish if either $X,Y\in H_\pm$
or $X\in H_\pm,Y\in Q_\pm$. If we insert $$ T^\pm(X,Y)=
\Theta^\pm(X,Y)-\Theta^\pm(Y,X),$$
in (\ref{equivCRFK}), we get (\ref{equivCRFK2}).
\end{proof}

The last form of the CRFK-conditions that we will prove is
\begin{prop}\label{CRFKcuLC} The generalized metric CRF-structure
$(\Phi,G)$ with the associated structures $(F_\pm,\gamma,\psi)$ is
a CRFK-structure iff $F_\pm$ are classical metric CRF-structures
and the following conditions hold for any $X,Y,Z\in\chi^1(M)$:
\begin{equation}\label{CRFK9}
\gamma(F_\pm X,(\nabla_ZF_\pm)(Y))= \pm\frac{1}{2}[d\psi(F^2_\pm
X,Y,Z)+d\psi(F_\pm X,F_\pm Y,Z)].\end{equation}
\end{prop}
\begin{proof}
We replace (\ref{equivCRFK2}) by conditions with general arguments
$X,Y,Z\in\chi^1(M)$ by replacing the arguments in $H_\pm$ by
arguments of the form $(F_\pm^2+\sqrt{-1}F_\pm)X$ and the
arguments in $Q_\pm$ by $(Id+F_\pm^2)Y$. After this replacements,
the conditions present a real and an imaginary part, which are
equivalent via the change $X\mapsto F_\pm X$, and we get the
following characteristic conditions of the CRFK-structures
\begin{equation}\label{CRFK5} \begin{array}{l}
\gamma(\Theta^\pm(Z,F_\pm Y), F_\pm^2X) +
\gamma(\Theta^\pm(Z,F_\pm^2 Y), F_\pm X)\vspace{2mm}\\
=\mp\frac{1}{2}[d\psi (F_\pm X,F_\pm^2Y,Z) +d\psi (F_\pm^2 X,F_\pm
Y,Z)],\vspace{2mm}\\
\gamma(\Theta^\pm(Z,Y+F_\pm^2 Y), F_\pm X) =
\mp\frac{1}{2}d\psi (F_\pm X,Y+F_\pm^2Y,Z).
\end{array}\end{equation}

Now, if we use (\ref{condTheta2}), (\ref{condTheta3}) and the
equality $$\nabla F_\pm^2=F_\pm\circ\nabla F_\pm+\nabla F_\pm\circ
F_\pm,$$ and then subtract the first condition from the second, we
get the following system that is equivalent to (\ref{CRFK5}):
\begin{equation}\label{CRFK6}
\begin{array}{l}
\gamma(F_\pm X,(\nabla_ZF_\pm)(F_\pm Y))=
\pm\frac{1}{2}[d\psi(F_\pm X,F_\pm^2Y,Z)+d\psi(F^2_\pm X,F_\pm
Y,Z)],\vspace{2mm}\\ \gamma(F_\pm X,F_\pm(\nabla_ZF_\pm)(Y)) =
\pm\frac{1}{2}[d\psi(F_\pm
X,Y,Z)-d\psi(F_\pm^2X,F_\pm Y,Z)],\end{array}\end{equation} which
does not contain the difference tensor $\Theta$ anymore. Moreover,
by replacing $X\mapsto F_\pm X, Y\mapsto F_\pm Y$ in the second
condition (\ref{CRFK6}), it follows that the second condition
implies the first. Finally, by means of the change $X\mapsto F_\pm
X$ we see that the (only remaining) second condition (\ref{CRFK6})
is equivalent with the required condition (\ref{CRFK9})
\end{proof}
\begin{rem}\label{obsXP} {\rm
In the generalized K\"ahler case, $F_\pm=J_\pm$ are complex
structures and the following formula, where $\omega_\pm$ are the
K\"ahler forms of the Hermitian structures $(\gamma,J_\pm)$, holds
\begin{equation}\label{formulaKN} \gamma((\nabla_XJ_\pm)(Y),Z) =
\frac{1}{2}[d\omega_\pm(X,Y,Z)-d\omega_\pm(X,J_\pm Y,J\pm
Z)]\end{equation} (e.g., Proposition IX.4.2 of
\cite{KN}; in (\ref{formulaKN}),
the constant factors are different from those in \cite{KN} because we
use different factor conventions for the exterior differential and
product and a different sign of the fundamental form).
Accordingly, (\ref{CRFK9}) becomes
\begin{equation}\label{auxptKGalt} \begin{array}{l}d\psi(X,Y,Z)-d\psi(J_\pm
X,J_\pm Y,Z)\vspace{2mm}\\=\pm[d\omega_\pm(J_\pm X,
Y,Z)+d\omega_\pm(X,J_\pm Y,Z)].\end{array}\end{equation} In
(\ref{auxptKGalt}), if we replace $X\mapsto J_\pm X$, then
subtract the first cyclic permutation of $(X,Y,Z)$ and add the
second cyclic permutation, we get
\begin{equation}\label{KG2}\begin{array}{l}
d\psi(X,Y,Z)=\pm\frac{1}{2}[d\omega_\pm(JX,JY,JZ)\vspace{2mm}\\ +
d\omega_\pm(JX,Y,Z)+d\omega_\pm(X,JY,Z)+d\omega_\pm(X,Y,JZ)].\end{array}
\end{equation}
Finally, by using arguments in the eigenbundles of $J_\pm$ and
since $\omega_\pm$ is of the complex type $(1,1)$ and
$d\omega_\pm$ has no $(3,0)$ and $(0,3)$ type components, it is
easy to check that the last three terms of (\ref{KG2}) add up to
the first term. Therefore, we get Gualtieri's characteristic
conditions for generalized K\"ahler structures
\cite{Galt}
\begin{equation}\label{eqluiGaltieri} d\psi(X,Y,Z)=\pm
d\omega_\pm(J_\pm X,J_\pm Y,J_\pm Z).
\end{equation}}\end{rem}
\begin{prop}\label{cazulpsi0} A generalized CRFK-manifold with a
closed $2$-form $\psi$ is a triple $(M,\gamma,\psi)$ where $\psi$
is a closed $2$-form and $\gamma$ is a Riemannian metric that has
two partially K\"ahler reductions.\end{prop}
\begin{proof} A Riemannian metric $\gamma$ is said to have a
partially K\"ahler reduction if there exists an atlas of local
coordinates $(z^a,y^u)$, where $z^a$ are complex and $y^u$ are
real, which has smooth, local transition functions
$$ \tilde z^a=\tilde z^a(z),\,\tilde y^u=\tilde y^u(y)$$ with
complex analytic functions $\tilde z^a(z)$, and
\begin{equation}\label{reduct}
\gamma=\gamma_{a\bar b}(z)dz^ad\bar
z^b+\gamma_{uv}(y)dy^udy^v\;\;(\gamma_{b\bar
a}=\bar{\gamma}_{a\bar b}),\end{equation} where the first term is
a K\"ahler metric.

If $d\psi=0$, the CRFK-conditions (\ref{CRFK9}) reduce to the
condition
\begin{equation}\label{CRFK10} \gamma(F_\pm X,(\nabla_ZF_\pm)(Y))=0,
\end{equation} which is equivalent with $(\nabla_ZF_\pm)(Y)\in Q_\pm$ for
any $Y,Z\in\chi^1(M)$.

If $Y\in Q_\pm$ the previous condition gives $F_\pm(\nabla_ZY)\in
Q_\pm\cap(im\,F_\pm)$, therefore, $F_\pm\nabla_Z(Y)=0$ and
$\nabla_ZY\in Q_\pm$. Thus, the distribution $Q_\pm$ is
$\nabla$-parallel, the same must hold for its $\gamma$-orthogonal
distribution $P_\pm=im\,F_\pm$, the distributions $P_\pm,Q_\pm$
are foliations and $\gamma$ is a reducible metric in two ways.

Furthermore, we get
$$(\nabla_ZF_\pm)(Y)=\nabla_Z(F_\pm Y)-F_\pm(\nabla_ZY)\in P_\pm,\;\;
\forall Y\in P_\pm,Z\in\chi^1(M),$$
hence, $\gamma(F_\pm X,(\nabla F_\pm)Y)=0$ implies
\begin{equation}\label{aux10}
(\nabla_ZF_\pm)(Y)=0\end{equation} for all $Y\in P_\pm$. By
looking at the case $Z\in Q_\pm$ we can see that the structures
$F_\pm|_{P_\pm}$ are projectable onto the space of leaves of the
foliation $Q_\pm$. Indeed, projectability means that for any
projectable vector field (an infinitesimal automorphism of the
foliation) $Y\in P_\pm$ the field $F_\pm Y$ is projectable too;
this holds since, $\forall Z\in\Gamma Q_\pm$ we have
$$[Z,F_\pm Y]=\nabla_Z(F_\pm Y)-\nabla_{F_\pm Y}Z
\stackrel{(\ref{aux10})}{=}F_\pm(\nabla_ZY)-\nabla_{F_\pm Y}Z$$
$$=F_\pm(\nabla_ZY-\nabla_YZ)-\nabla_{F_\pm Y}Z=F[Z,Y] -\nabla_{F_\pm Y}Z
=-\nabla_{F_\pm Y}Z\in Q_\pm.$$ (If $Z\in Q_\pm$ then
$\nabla_YZ\in Q_\pm$ and $F\nabla_YZ=0$. On the other hand, the
projectability of $Y$ means that $[Z,Y]\in Q_\pm$ and
$F_\pm[Z,Y]=0$.) Then, condition (\ref{aux10}) for $Z\in P_\pm$
exactly means that $(\gamma|_{P_\pm},F_\pm|_{P\pm})$ is a K\"ahler
structure on the leaves of foliation $P_\pm$, hence, the
reductions of the metric $\gamma$ mentioned above are partially
K\"ahler reductions. The chain of arguments may be reversed. This
leads from the partial K\"ahler reductions and the closed form
$\psi$ to a CRFK-structure. \end{proof}
\begin{example}\label{exfinal1} {\rm
If $(F,\gamma)$ is a classical metric F-structure and $(\Phi,G)$
the corresponding generalized structure given in Example
\ref{Fmetricclasic}, the latter is CRFK iff the classical
structure $(F,\gamma)$ comes from a partial K\"ahler manifold.
Accordingly, the generalized CRFK-manifolds should be seen as
generalized, partially K\"ahler manifolds.}\end{example}
\begin{example}\label{exfinal2} {\rm Take $M=\mathds{R}^{2n+h}$
with the canonical coordinates $(x^i,y^u)$
$(i=1,...,2n;u=1,...,h)$ and with the Euclidean metric
$$\gamma=\sum_{i=1}^{2n}(dx^i)^2 + \sum_{u=1}^{h}(dy^u)^2.$$
If we define complex coordinates $z^a=x^a+\sqrt{-1}x^{n+a},\;\;
a=1,...,n$, $\gamma$ gets the partially K\"ahler reduction
$$\gamma=\sum_{a=1}^{n}dz^adz^{\bar b} + \sum_{u=1}^{h}(dy^u)^2.$$
On the other hand, if we define complex coordinates
$w^u=x^u+\sqrt{-1}y^u$, we get a second partially K\"ahler
reduction
$$\gamma=\sum_{a=1}^{n}dw^udw^{\bar v} + \sum_{\alpha=h+1}^{2n}(dx^\alpha)^2.$$
Therefore, for any closed $2$-form $\psi$, $(M,\gamma,\psi)$ is
endowed with a generalized CRFK-structure. By the usual
quotientizing, the previous structure a generalized CRFK-structure
on the torus
$\mathds{T}^{2n+h}=\mathds{R}^{2n+h}/\mathds{Z}^{2n+h}$.}\end{example}
\section{Appendix: Generalized Sasakian Structures} The structures
discussed in Section 4 do not include properly generalized Sasakian
manifolds \cite{Bl1}. In view of the importance of the latter we
define such a generalization in this Appendix. With the notation of
Example \ref{exgencontact}, consider a generalized almost contact
structure of codimension $h=1$. This structure may be identified
with the generalized almost complex structure of $M\times\mathds{R}$
that has the classical tensor fields
\begin{equation}\label{structura'} A'= F,\,
\pi'=P+Z\wedge\frac{\partial}{\partial t},\,\sigma'=\theta+\xi\wedge
dt.\end{equation} More exactly, it is easy to see that the
generalized almost complex structures $\Phi$ of $M\times\mathds{R}$
that can be obtained in this way are those that satisfy the
following two properties:

(i) $\Phi$ is invariant by translations along $\mathds{R}$,

(ii) $\Phi(T\mathds{R}\oplus0)\subseteq 0\oplus
T^*M,\;\Phi(0\oplus T^*\mathds{R})\subseteq TM\oplus0$.\\

\noindent The addition of the property

(iii) $\Phi(TM\oplus0)\subseteq TM\oplus T^*\mathds{R},\;
\Phi(0\oplus T^*M)\subseteq T\mathds{R}\oplus T^*M,$\\

\noindent takes us to the case $P=0,\theta=0$, which is the case of
a classical almost contact structure. This is an interpretation of
the classical almost contact structures of $M$ by non-classical,
generalized, almost complex structures of $M\times\mathds{R}$.

In classical geometry the almost contact structure $(F,Z,\xi)$ is
identified with the classical almost complex structure of
$M\times\mathds{R}$ defined by \cite{Bl1}
\begin{equation}\label{Jacontact} J=F+dt\otimes
Z-\xi\otimes\frac{\partial}{\partial t}.\end{equation} The
structures $J$ obtained in this way are characterized by the
properties

(i) $J$ is translation invariant,

(ii) $J(T\mathds{R})\subseteq TM$.\\

We also notice the important fact that the integrability of $J$ is
equivalent with the normality of $(F,Z,\xi)$ and we know from
Example \ref{exgencontact} that normality is also equivalent with
the integrability of the corresponding generalized almost complex
structure $\Phi_0$ given by (\ref{structura'}) with
$P=0,\theta=0$.

The addition of a Riemannian metric $\gamma$ of $M$ such that
\begin{equation}\label{acmetric} \gamma(FX,FY)=\gamma(X,Y)-\xi(X)\xi(Y)
\end{equation}
(which implies $\xi=\flat_\gamma Z,g(Z,FX)=0,g(Z,Z)=1$ \cite{Bl1})
yields an almost contact metric structure $(F,Z,\xi,\gamma)$. It is
easy to check that (\ref{acmetric}) is equivalent with the fact that
the metric
\begin{equation}\label{gamma}
\Gamma=e^{t}(\gamma+dt^2)\end{equation} is Hermitian for the
almost complex structure $J$ of $M\times\mathds{R}$ defined by
(\ref{Jacontact}). The factor $e^{t}$ is superfluous here but
essential for the notion of a Sasakian structure. Thus, the almost
contact metric structures of $M$ may be seen as translation
invariant, almost Hermitian structures $(\Gamma,J)$ of
$M\times\mathds{R}$ (property (ii) of $J$ necessarily holds since
$J(\partial/\partial t)\perp_{\Gamma}\partial/\partial t$).

The almost contact metric structure $(F,Z,\xi,\gamma)$ has the
associated fundamental $2$-form $\Xi(X,Y)=g(FX,Y)$ while the
corresponding almost Hermitian structure $J$ has the K\"ahler form
$\omega$. A simple calculation gives
\begin{equation}\label{formaomega}
\omega=e^{t}(\Xi-\xi\wedge
dt),\;d\omega=e^{t}[d\Xi+(\Xi-d\xi)\wedge dt].\end{equation}

The most usual definition of a Sasakian structure (e.g.,
\cite{Bl1}) requires it to be a normal, contact, metric structure
$(F,Z,\xi,\gamma)$ where the use of the term {\it contact} instead
of {\it almost contact} means the requirement $\Xi=d\xi$. Thus, in
view of (\ref{formaomega}), a Sasakian structure is characterized
by the fact that the corresponding structure $(\Gamma,J)$ is
K\"ahler
\cite{Bl1}.

The previous remark suggests the question of determining the
almost contact metric structures $(F,Z,\xi,\gamma)$ such that the
corresponding structure $\Phi_0$ is generalized K\"ahler for a
convenient metric, which turns out to be $\bar\Gamma=\gamma+dt^2$.
\begin{prop}\label{caractcosympl} The
almost contact metric structure $(F,Z,\xi,\gamma)$ corresponds to
a generalized K\"ahler structure $(\bar\Gamma,\Phi_0)$ iff the
structure is cosymplectic in the sense of Blair.\end{prop}
\begin{proof} Here, $\bar\Gamma$ is to be seen as the
generalized Riemannian metric of $M\times\mathds{R}$ given by
(\ref{matriceaG}) with $\varphi=0,\gamma=\beta=\bar\Gamma$. Then,
using Proposition \ref{metricaHermite}, it is easy to check that the
metric conditions (\ref{acmetric}) are equivalent with the fact that
$(\bar\Gamma,\Phi_0)$ is a generalized almost Hermitian structure.

The generalized almost Hermitian structure $(\bar\Gamma,\Phi_0)$
has two associated classical almost complex structures defined by
formula (\ref{formuleJ}) with $\psi=0$. These structures are
$$J_\pm=A'\pm\sharp_{\pi'}\circ\flat_{\bar\Gamma}=F\mp dt\otimes Z
\pm\xi\otimes\frac{\partial}{\partial t}$$ ($J_-$ is exactly the
structure (\ref{Jacontact})). The K\"ahler forms of
$(\bar\Gamma,J_\pm)$ are $$\omega_\pm=\Xi\pm\xi\wedge dt.$$ With
(\ref{eqluiGaltieri}), the structure $(\bar\Gamma,\Phi_0)$ is
generalized K\"ahler iff $J_\pm$ are integrable and
$d\omega_\pm=0$. Therefore, $(\bar\Gamma,\Phi_0)$ is generalized
K\"ahler iff the structure $(F,Z,\xi,g)$ is normal and satisfies
the conditions
\begin{equation}\label{cosympl} d\xi=0,\;d\Xi=0.\end{equation}
This exactly is the definition of a cosymplectic structure in the
sense of Blair \cite{Bl1}.\end{proof}

It is known that a cosymplectic structure in the sense of Blair may
be identified with a partially K\"ahler metric of the form
(\ref{reduct}) where the second term is $dt^2$ \cite{Bl1}.
Therefore, as noticed in Example \ref{exfinal1}, the generalized
CRFK-structures included adequate, generalized, Blair-cosymplectic
structures.

Earlier, we saw that a Sasakian manifold is a Riemannian manifold
$(M,\gamma)$ endowed with a translation invariant complex structure
$J$ of $M\times\mathds{R}$ such that $(M\times\mathds{R},\Gamma,J)$
is a K\"ahler manifold. Now, let us assume that $M$ is endowed with
a generalized Riemannian metric $G$, equivalently, with a pair
$(\gamma,\psi)$ where $\gamma$ is a classical Riemannian metric and
$\psi\in\Omega^2(M)$. Then, the generalized Riemannian metrics
$\tilde G$ of $M\times\mathds{R}$ that are related to $G$ in the way
$\Gamma$ was related to $\gamma$ are those defined by pairs
\begin{equation}\label{tildeG} \tilde G\Leftrightarrow
(\Gamma,\Psi=e^{t}(\psi+\kappa\wedge dt))\end{equation} where
$\kappa$ is an arbitrary $1$-form on $M$. Accordingly, we give
\begin{defin}\label{defgenS} {\rm A {\it generalized Sasakian manifold}
is a generalized Riemannian manifold $(M,(\gamma,\psi))$ endowed
with a translation invariant generalized complex structure $\Phi$ of
$M\times\mathds{R}$ such that, for some $\kappa\in\Omega^1(M)$,
$(M\times\mathds{R},\tilde G,\Phi)$ is a generalized K\"ahler
manifold.}\end{defin}

In this definition $\tilde G$ is defined by (\ref{tildeG}) and the
invariance of $\Phi$ by translations means $L_{\partial/\partial
t}\Phi=0$, where the Lie derivative is defined like for an
endomorphism of $TM$ and acts on both the $TM$ and $T^*M$
components. In the next proposition we will also use the following
notation: $\forall\lambda\in\Omega^k(M)$, $\lambda^c$ will be the
form obtained by evaluating $\lambda$ on arguments of the form
$F_\pm X$.
\begin{prop}\label{caractSgen} A generalized Sasakian structure of
a manifold $M$ is equivalent with a pair of classical, normal,
almost contact, metric structures $(F_\pm,Z_\pm,\xi_\pm,\gamma)$
complemented by a pair of forms
$\psi\in\Omega^2(M),\kappa\in\Omega^1(M)$ that satisfy  the
conditions
\begin{equation}\label{bGualt8}\begin{array}{c}
i(Z_\pm)(\psi+d\kappa)=0,\,(\psi+d\kappa)^c
=-L_{Z_\pm}[L_{Z_\pm}(\psi+d\kappa)]^c,\vspace{2mm}\\
(d\psi)^c=-i(Z_\pm)[\xi_\pm\wedge
d[L_{Z_\pm}(\psi+d\kappa)]^c],\end{array}\end{equation} and
\begin{equation}\label{bGualt9}
\Xi_\pm=d\xi_\pm\mp [L_{Z_\pm}(\psi+d\kappa)]^c.\end{equation}
\end{prop}
\begin{proof} The structure $\Phi$ of Definition \ref{defgenS} is
equivalent with a pair of $\Gamma$-compatible, translation
invariant, classical, complex structures $J_\pm$ of
$M\times\mathds{R}$ and we have seen that, in turn, this pair is
equivalent with a pair of normal, almost contact metric structures
$(F_\pm,Z_\pm,\xi_\pm,\gamma)$ of $M$. The remaining conditions to
be discussed are Gualtieri's conditions (\ref{eqluiGaltieri}) for
the K\"ahler forms of $(\Gamma,J_\pm)$ given by
(\ref{formaomega}).

These conditions are
\begin{equation}\label{bGualt} \begin{array}{l}d\omega_\pm(J_\pm
(X+a\frac{\partial}{\partial t}), J_\pm
(Y+b\frac{\partial}{\partial t}),J_\pm(U+u\frac{\partial}{\partial
t})\vspace{2mm}\\=\pm d\Psi(X+a\frac{\partial}{\partial
t},Y+b\frac{\partial}{\partial t},U+u\frac{\partial}{\partial
t}).\end{array}\end{equation} Using (\ref{Jacontact}) and
(\ref{formaomega}), (\ref{bGualt}) becomes
\begin{equation}\label{bGualt1}
\begin{array}{l} d\Xi_\pm(F_\pm X,F_\pm Y,F_\pm U)
+\sum_{Cycl}u[i(Z_\pm)d\Xi_\pm](F_\pm X,F_\pm Y)\vspace{2mm}\\
-\sum_{Cycl}\xi_\pm(U)(\Xi_\pm-d\xi_\pm)(F_\pm X+aZ_\pm,F_\pm Y
+bZ_\pm)\vspace{2mm}\\ =\pm\{d\psi(X,Y,U)+
\sum_{Cycl}u(\psi+d\kappa)(X,Y)\},
\end{array}\end{equation} where the cyclic permutations in the
sums are on the arguments $(X,a),(Y,b),(U,u)$.

Since $T(M\times\mathds{R})=TM\oplus T\mathds{R}$ and
(\ref{bGualt}) is the equality of $3$-forms, it follows that
(\ref{bGualt1}) holds iff it holds in the cases 1)
$a=0,b=0,u=1,U=0$ and 2) $a=b=u=0$. In case 1), (\ref{bGualt1})
reduces to
\begin{equation}\label{bGualt2}\begin{array}{c} [i(Z_\pm)d\Xi_\pm](F_\pm X,F_\pm
Y)+\{\xi_\pm\wedge[i(Z_\pm)(\Xi_\pm-d\xi_\pm)]\circ F_\pm\}(X,Y)
\vspace{2mm}\\
=\pm(\psi+d\kappa)(X,Y),\end{array}\end{equation} which, by taking
into account $i(Z_\pm)\Xi_\pm=0$ and the fact that
$\xi_\pm(Z_\pm)=1$ and normality imply
$i(Z_\pm)d\xi_\pm=L_{Z_\pm}\xi_\pm=0$ \cite{Bl1}, becomes
\begin{equation}\label{bGualt20}[L_{Z_\pm}\Xi_\pm](F_\pm X,F_\pm
Y)=\pm(\psi+d\kappa)(X,Y).\end{equation} In case 2),
(\ref{bGualt1}) reduces to
\begin{equation}\label{bGualt3} \begin{array}{c}
d\Xi_\pm(F_\pm X,F_\pm Y,F_\pm U)
-\sum_{Cycl}\xi_\pm(U)[(\Xi_\pm-d\xi_\pm)(F_\pm X,F_\pm
Y)]\vspace{2mm}\\ =\pm d\psi(X,Y,U).\end{array}\end{equation}

Furthermore, since $TM=im\,F_\pm\oplus span\{Z_\pm\}$ we may
decompose each condition (\ref{bGualt20}), (\ref{bGualt3}) into
the cases (i) one of the arguments is $Z_\pm$ and the others
belong to $im\,F_\pm$, (ii) all the arguments are in $im\,F_\pm$.
In case (i) we have to write (\ref{bGualt20}) for
$(X,Y)\mapsto(F_\pm X,Z_\pm)$ and the result is
\begin{equation}\label{bGualt4}[i(Z_\pm)(\psi+d\kappa)]^c=0.
\end{equation} Since $<i(Z_\pm)(\psi+d\kappa),Z_\pm>=0$,
(\ref{bGualt4}) is the first condition (\ref{bGualt8}). In case
(ii) we have to write (\ref{bGualt20}) for $(X,Y)\mapsto(F_\pm
X,F_\pm Y)$ and the result is
\begin{equation}\label{bGualt5}
L_{Z_\pm}\Xi_\pm=\pm(\psi+d\kappa)^c.
\end{equation} Similarly, we have to express
(\ref{bGualt3}) for the cases (i) $(X,Y,U)\mapsto (F_\pm X,F_\pm
Y,Z_\pm)$ (ii) $(X,Y,U)\mapsto (F_\pm X,F_\pm Y,F_\pm U)$. For (i)
we get
\begin{equation}\label{bGualt6}
\Xi_\pm-d\xi_\pm=\mp [i(Z_\pm)d\psi]^c\end{equation} and
for (ii) we get
\begin{equation}\label{bGualt7}
d\Xi_\pm-\xi_\pm\wedge(i(Z_\pm)d\Xi_\pm) =i(Z_\pm)(\xi_\pm\wedge
d\Xi_\pm) =\mp (d\psi)^c,\end{equation} where the first equality
follows from the properties of the operator $i(Z_\pm)$.

Furthermore, we may calculate $L_{Z_\pm}\Xi_\pm$ and $d\Xi_\pm$
from (\ref{bGualt6}) and insert the corresponding values in
(\ref{bGualt5}), (\ref{bGualt7}). The result will be a system of
conditions that look exactly like (\ref{bGualt8}),
(\ref{bGualt9})except for the fact that instead of the form
$L_{Z_\pm}(\psi+d\kappa)$ one has $i(Z_\pm)d\psi$. But, modulo
(\ref{bGualt4}), we may replace the former by the latter.
\end{proof}
\begin{corol}\label{concluzie2} If $\psi$ is closed the structures
$(F_\pm,Z_\pm,\xi_\pm,\gamma)$ of a generalized Sasakian manifold
$M$ are classical Sasakian structures.\end{corol}
\begin{proof} Use $L_{Z_\pm}(\psi+d\kappa)=i(Z_\pm)d\psi$ in
(\ref{bGualt8}), (\ref{bGualt9}).\end{proof}
\begin{corol}\label{concluzie1} Let $(F_\pm,Z_\pm,\xi_\pm,\gamma)$ be a
pair of normal almost contact metric structures and let
$\psi\in\Omega^2(M),\kappa\in\Omega^1(M)$ be forms such that
\begin{equation}\label{cazbazic} i(Z_\pm)(\psi+d\kappa)=0,\;
L_{Z_\pm}(\psi+d\kappa)=0. \end{equation} Then, these data define
a generalized Sasakian structure on $M$ iff the structures
$(F_\pm,Z_\pm,\xi_\pm,\gamma)$ are classical Sasakian structures.
\end{corol}
\begin{proof} Insert (\ref{cazbazic}) in (\ref{bGualt8}), (\ref{bGualt9}).
\end{proof}

We should notice that a generalized Sasakian structure in the
sense of Definition \ref{defgenS} may not be a generalized almost
contact structure as defined in the first part of this Appendix.
Recall that the latter was characterized by the properties (i),
(ii) mentioned after formula (\ref{structura'}) and if we use
formula (\ref{eqpiAdinF}) for the structures $J_\pm$ of a
generalized Sasakian manifold (instead of $F_\pm$) to reconstruct
the generalized complex structure $\Phi$ we see that $\Phi$ may
not satisfy the required property (ii). However, $\Phi$ will
satisfy property (ii) if $\kappa=0,Z_-=-Z_+,\xi_-=-\xi_+$. A
similar terminological difficulty is encountered if one defines
the notion of a generalized, metric, almost contact structure by a
translation invariant, generalized, almost complex structure
$\Phi$ of $M\times\mathds{R}$ that is Hermitian with respect to
the generalized metric $\tilde G$, equivalently, by a pair of
classical metric, almost contact structures $(F_+,Z_+,\xi_+)$,
$(F_-,Z_-,\xi_-)$ with the same metric $\gamma$ and a pair of
forms $(\psi\in\Omega^2(M),\kappa\in\Omega^1(M))$.
\hspace*{7.5cm}{\small \begin{tabular}{l} Department of
Mathematics\\ University of Haifa, Israel\\ E-mail:
vaisman@math.haifa.ac.il \end{tabular}}
\end{document}